\documentclass[]{amsart}

\usepackage{bbding}
\usepackage{graphicx}
\usepackage{stackrel,amsfonts,amssymb,amsmath,amsthm,amscd}
\allowdisplaybreaks
\usepackage{caption}
\usepackage{subcaption}
\usepackage{mathrsfs}
\usepackage{xcolor}
\usepackage{empheq}
\usepackage{adforn}
\usepackage{accents}
\usepackage{cancel}
\usepackage{xr}
\usepackage{mdframed}
\usepackage{tikz-cd}
\usepackage{tikz}
\usepackage{amsfonts}
\usepackage{mathdots}
\usepackage{pinlabel}

\usepackage{enumerate}
\usepackage{lmodern}
\usepackage{mdframed}
\usepackage{stmaryrd}
\usepackage{blindtext}
\usepackage{leftidx}
\usepackage{easytable}
\usepackage{enumitem}
\usepackage{marginnote} 
\usepackage{easyReview}

\renewcommand{\thefootnote}{\ifcase\value{footnote}\or$\boldsymbol{\dagger}$\or
	$\boldsymbol{\ddagger}$\or$\boldsymbol{\dagger\dagger}$\or$\boldsymbol{\dagger\ddagger}$or$\boldsymbol{\ddagger\ddagger}$\fi}
\usepackage{perpage}
\MakePerPage{footnote}

\usepackage{tgpagella}
\usepackage[T1]{fontenc}

\usepackage[marginparwidth=2.5cm]{geometry}
\usepackage[backref=page]{hyperref}
\def\equationautorefname~#1\null{(#1)\null}
\usepackage{xcolor}
\hypersetup{
	colorlinks,
	linkcolor={[RGB]{0,75,150}},
	citecolor={[RGB]{160,0,0}},
	urlcolor={[RGB]{160,0,0}}
}
\PassOptionsToPackage{linktocpage}{hyperref}

\newcommand{\rleftarrows}{\mathrel{\raise.75ex\hbox{\oalign{%
				$\hfil\scriptstyle\relbar$\cr
				\vrule width0pt height.5ex$\scriptstyle\smash\leftarrow$\cr}}}}
\newcommand{\rightlarrows}{\mathrel{\raise.75ex\hbox{\oalign{%
				$\scriptstyle\rightarrow$\hfil\cr
				$\scriptstyle\vrule width0pt height.5ex\relbar$\cr}}}}
\newcommand{\Rrelbar}{\mathrel{\raise.75ex\hbox{\oalign{%
				$\scriptstyle\relbar$\cr
				\vrule width0pt height.5ex$\scriptstyle\relbar$}}}}

\makeatletter
\def\rightleftarrowsfill@{\arrowfill@\rleftarrows\Rrelbar\rightlarrows}
\newcommand{\xrightleftarrows}[2][]{\ext@arrow 3399\rightleftarrowsfill@{#1}{#2}}
\makeatother

\makeatletter
\newtheorem*{rep@theorem}{\rep@title}
\newcommand{\newreptheorem}[2]{%
	\newenvironment{rep#1}[1]{%
		\def\rep@title{#2 \ref{##1}}%
		\begin{rep@theorem}}%
		{\end{rep@theorem}}}
\makeatother

\makeatletter
\newtheorem*{rep@cor}{\rep@title}
\newcommand{\newrepcor}[2]{%
	\newenvironment{rep#1}[1]{%
		\def\rep@title{#2 \ref{##1}}%
		\begin{rep@cor}}%
		{\end{rep@cor}}}
\makeatother

\makeatletter
\newtheorem*{rep@prop}{\rep@title}
\newcommand{\newrepprop}[2]{%
	\newenvironment{rep#1}[1]{%
		\def\rep@title{#2 \ref{##1}}%
		\begin{rep@prop}}%
		{\end{rep@prop}}}
\makeatother
\newtheorem{corollary}{Corollary}[section]
\newtheorem{corx}{Corollary}

\newtheorem{theorem}[corollary]{Theorem}
\newtheorem{thmx}[corx]{Theorem}

\newtheorem{proposition}[corollary]{Proposition}

\newrepcor{cor}{Corollary}
\newreptheorem{theorem}{Theorem}
\newrepprop{prop}{Proposition}

\newtheorem*{theorem*}{Theorem}
\newtheorem{lemma}[corollary]{Lemma}

\theoremstyle{definition}

\newtheorem*{remark*}{Remark}
\newtheorem{definition}[corollary]{Definition}
\newtheorem*{definition*}{Definition}
\newtheorem{example}[corollary]{Example} 
\newtheorem{remark}[corollary]{Remark}

\theoremstyle{remark} 
\numberwithin{equation}{section}
\numberwithin{figure}{section}

\newcommand{\wh}{\hat}
\newcommand{\wt}{\widetilde}

\newcommand{\id}{\mathrm{id}}

\newcommand{\Hom}{\mathrm{Hom}}

\newcommand{\M}{\mathfrak{M}}

\newcommand{\ve}[1]{{\boldsymbol{#1}}}
\newcommand{\T}{\mathcal{T}}
\newcommand{\R}{\mathbb{R}}

\newcommand{\RP}{\mathbb{RP}}
\renewcommand{\P}{\mathbb{P}}
\newcommand{\C}{\mathbb{C}}
\newcommand{\CP}{\mathbb{CP}}
\renewcommand{\H}{\mathbb{H}}
\newcommand{\E}{\mathbb{E}}
\renewcommand{\SS}{\mathbb{S}}
\newcommand{\TT}{\mathscr{T}}
\newcommand{\dT}{\widetilde{\mathscr{T}}}

\newcommand{\pro}{\raisebox{-0.06cm}{\scalebox{1.6}{$\cdot$}}}

\newcommand{\Face}{\mathsf{F}}
\newcommand{\Edge}{\mathsf{E}}

\newcommand{\EP}{E}

\newcommand{\dif}{\mathsf{d}}

\newcommand{\transp}[1]{\leftidx{^\mathsf{t}}{#1}}

\newcommand{\pa}{\partial}

\newcommand{\dx}{{\dif x}}

\newcommand{\Isom}{\mathsf{Isom}}

\newcommand{\PSL}{\mathrm{PSL}}

\newcommand{\SO}{\mathrm{SO}}

\begin{document}

\title{Boundary metric of Epstein-Penner convex hull and discrete conformality} 

\author{Xin Nie}
\address{
	Shing-Tung Yau Center of Southeast University, Nanjing 210096, China}
\email{nie.hsin@gmail.com}

\maketitle

\begin{abstract}
The Epstein-Penner convex hull construction associates to every decorated punctured hyperbolic surface a polyhedral convex body in the Minkowski space. It works in the de Sitter and anti-de Sitter spaces as well. In these three spaces, the quotient of the spacelike boundary part of the convex body has an induced Euclidean, spherical and hyperbolic metric, respectively, with conical singularities. We show that this gives a bijection from the decorated Teichm\"uller space to a moduli space of such metrics in the Euclidean and hyperbolic cases, as well as a bijection between specific subspaces of them in the spherical case. Moreover, varying the decoration of a fixed hyperbolic surface corresponds to a discrete conformal change of the metric. This gives a new $3$-dimensional interpretation of discrete conformality which is in a sense inverse to the Bobenko-Pinkall-Springborn interpretation.
\end{abstract}

\tableofcontents

\section{Introduction}\label{sec_1}
Let $S_{g,n}$ be the closed orientable topological surface of genus $g$ with $n\geq1$ marked points, $S^\circ_{g,n}$ be the punctured surface obtained by removing the marked points, and $\TT_{g,n}$ be the  Teichm\"uller space formed by equivalence classes 
of finite-area complete hyperbolic metrics\footnote{
	The \emph{equivalence} of metrics on $S_{g,n}$ or $S^\circ_{g,n}$ is defined by isotopies fixing the marked points. 
	Adopting a common abuse of notation, we often do not distinguish a metric and the equivalence classe represented by it. For example, by writing ``$\ve{h}\in\TT_{g,n}$'', we mean that $\ve{h}$ is a finite-area complete hyperbolic metric on $S^\circ_{g,n}$, and similarly for other spaces of metrics introduced below.} on $S_{g,n}^\circ$.
The \emph{decorated Teichm\"uller space} $\dT_{g,n}$ is formed by pairs $(\ve{h},c)$, where $\ve{h}\in\TT_{g,n}$ is a hyperbolic metric and $c=(c_1,\cdots,c_n)$ is a choice of horocycles (a ``decoration'') centered at each marked point (see \cite{Penner}). The natural projection
$$
\Pi:\dT_{g,n}\to\TT_{g,n},\quad (\ve{h},c)\mapsto\ve{h}
$$
gives $\dT_{g,n}$ the structure of a fiber bundle over $\TT_{g,n}$, where each fiber $\Pi^{-1}(\ve{h})$ is a product of $n$ affine lines. 

\begin{figure}[h]
 	\centering
 	\includegraphics[width=12.5cm]{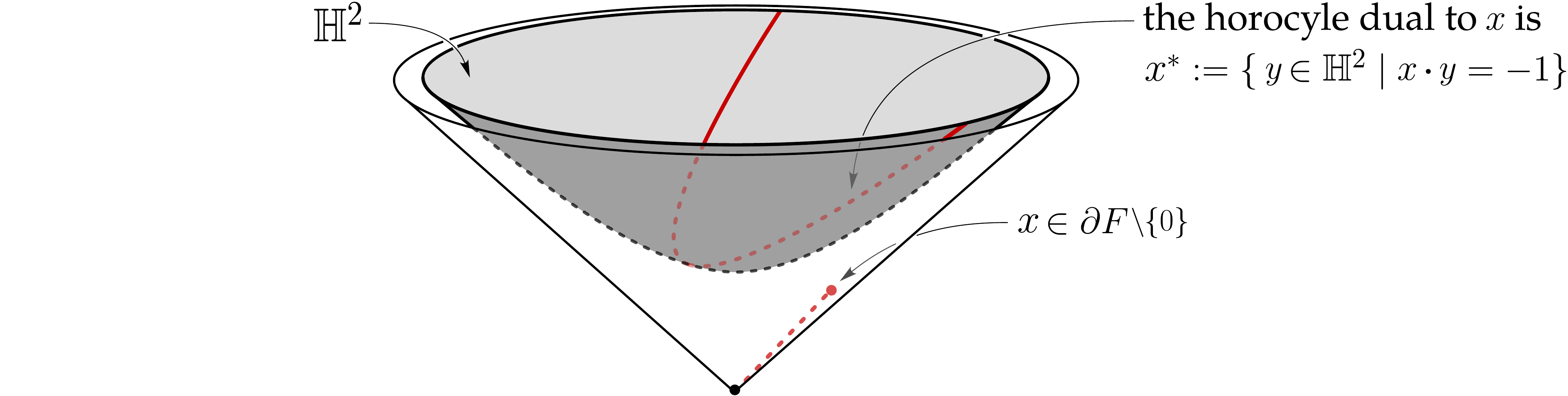}
 	\caption{Duality between points of $\pa F\setminus\{0\}$ and horocycles in $\H^2$.}
 	\label{figure_duality}
 \end{figure}

\begin{figure}[h]
	\centering
	\includegraphics[width=8.1cm]{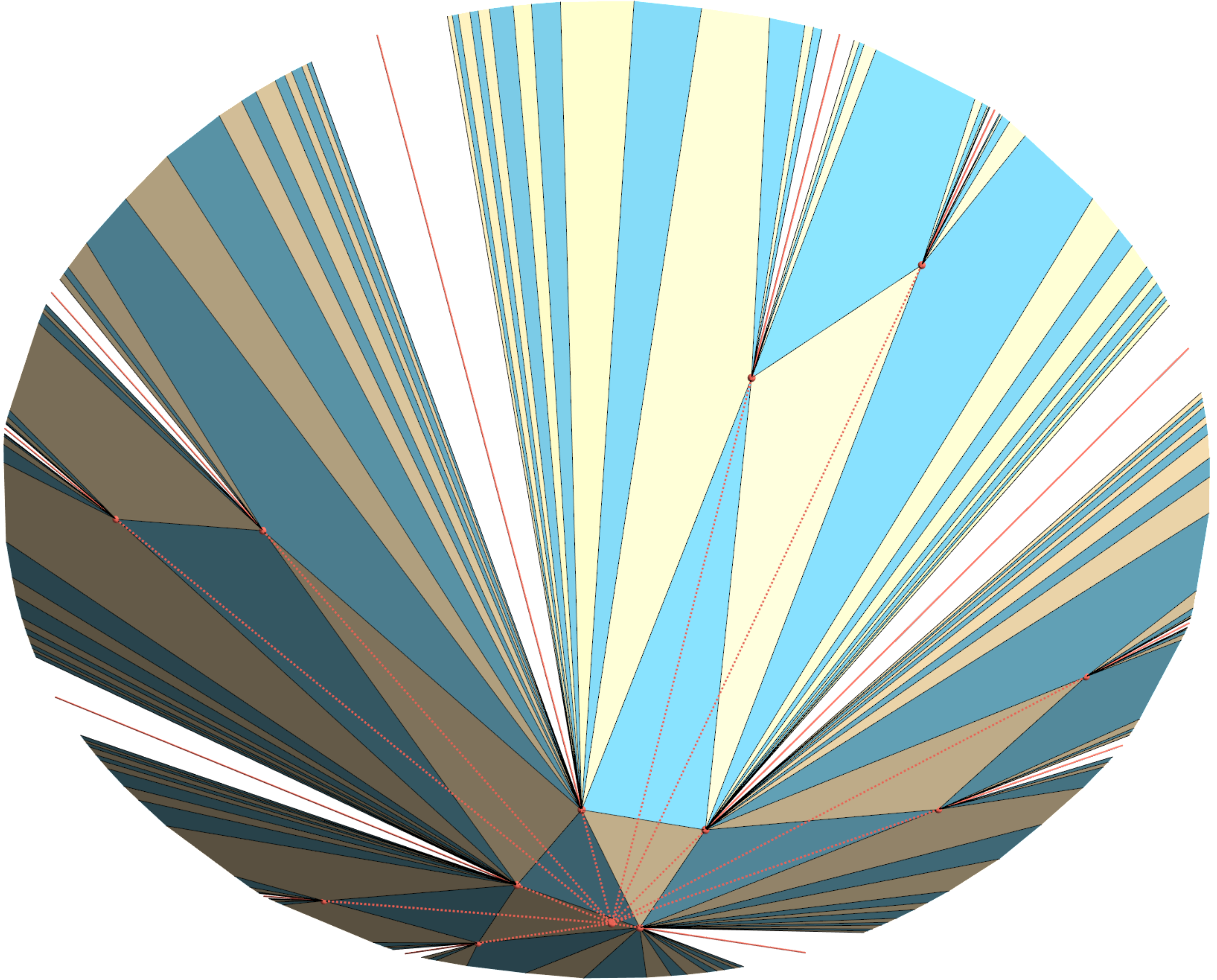}
	\hspace{1.8cm}
	\includegraphics[width=6.2cm]{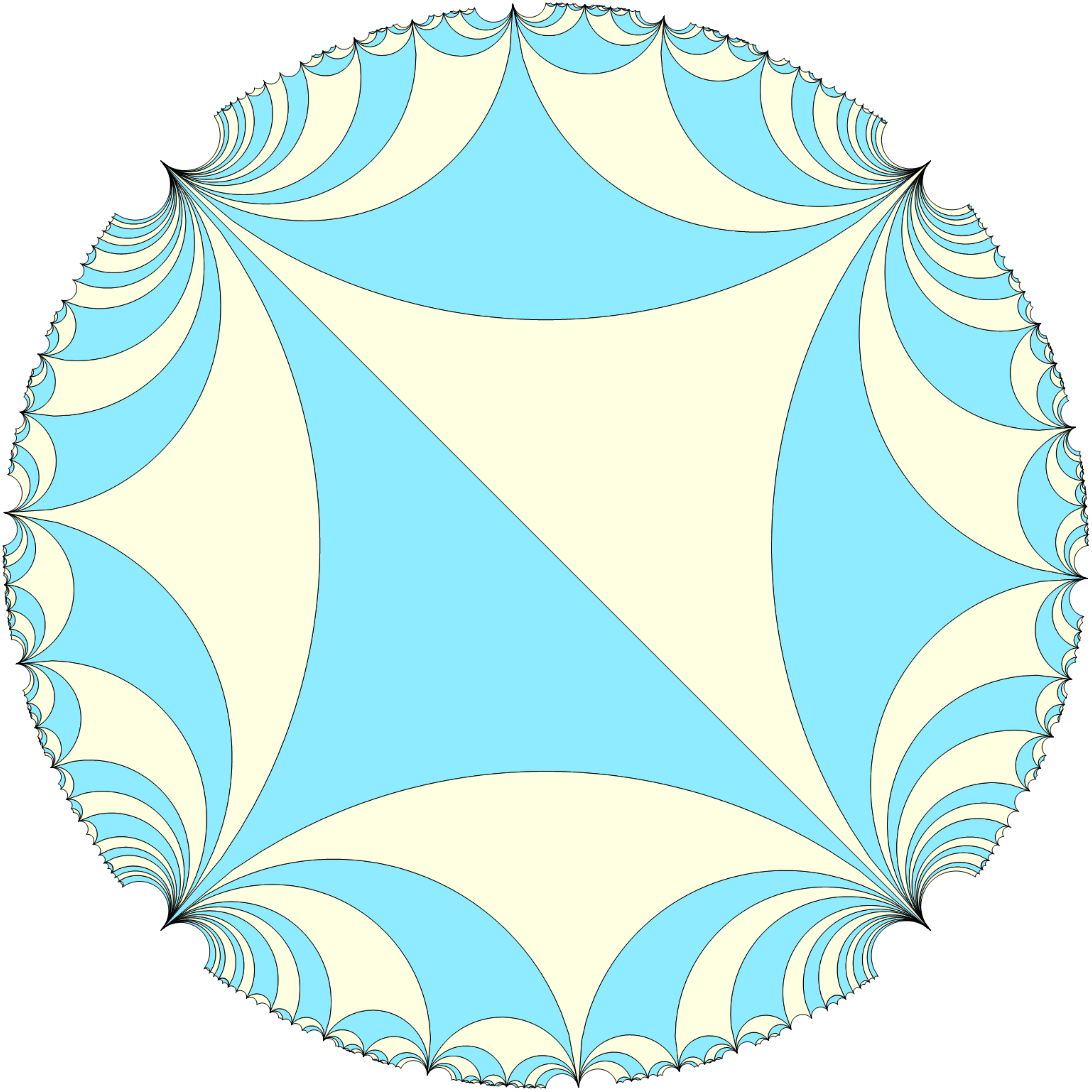}
	\caption{Boundary of an Epstein-Penner convex hull, viewed from below, with adjacent faces colored differently. There are infinitely many faces at each vertex, accumulating along a lightlike ray, while the figure only shows finitely many of them. On the right is its image in $\H^2$ (identified with the Poincar\'e disk) by radial projection. See also Figure \ref{figure_halfpipe}.
	}
	\label{figure_hull}
\end{figure}

Let $F:=\{x\in\R^{2,1}\mid x\pro x\leq0,\,x_3\geq0\}$ be the future cone of the origin in the Minkowski space $\R^{2,1}$. Given a decorated metric $(\ve{h},c)\in\dT_{g,n}$, the \emph{Epstein-Penner convex hull} $\EP(\ve{h},c)$ is a convex body in $F$, well defined up to $\SO_0(2,1)$ (the Minkowski motions preserving $F$), constructed as follows (\cite{Epstein-Penner, Penner}): 
Identify $(S_{g,n}^\circ,\ve{h})$ with a quotient of the hyperboloid $\H^2=\{x\pro x=-1,x_3>0\}$ by a Fuchsian group $\Gamma\subset\SO_0(2,1)$.
 Under the one-to-one correspondence between points of $\pa F\setminus\{0\}$ and horocycles in $\H^2$ as shown in Figure \ref{figure_duality}, 
 there is a $\Gamma$-invariant set $C\subset\pa F\setminus\{0\}$ whose corresponding horocycles are the lifts of the  decorating horocycles $c$. Then $\EP(\ve{h},c)$ is defined to be the convex hull of $C$ in $\R^{2,1}$, see Figure \ref{figure_hull}. 
It is shown in \cite{Epstein-Penner} that the boundary $\pa \EP(\ve{h},c)$ consists of the following two parts, which overlap only at $C$:
\begin{itemize}
	\item the lightlike part $\pa \EP(\ve{h},c)\cap\pa F$ is the union of rays $\{\alpha x\mid \alpha\geq1,\, x\in C\}$;
	\item the spacelike part $\pa^+\! \EP(\ve{h},c):=(\pa \EP(\ve{h},c)\setminus\pa F)\cup C$ is a union of spacelike planar polygons with vertices in $C$. For a generic $(\ve{h},c)$, these polygons are all triangles.
\end{itemize}
The quotient of $\pa^+\! \EP(\ve{h},c)$ by $\Gamma$ naturally identifies with $S_{g,n}$, so the induced metric on $\pa^+\! \EP(\ve{h},c)$ descends to a Euclidean metric on $S_{g,n}$ which may have conical singularities at the marked points. This yields a map from $\dT_{g,n}$ to the space $\M^0_{g,n}$ of equivalence classes of such metrics.


More generally, let $\M^K_{g,n}$ be the space of equivalence classes of metrics on $S_{g,n}$ of curvature $K\in\{-1,0,1\}$, possibly with conical singularities at marked points.
A new observation of this paper is that the Epstein-Penner construction is still valid when $F$ is replaced by the future cone of any point in the \emph{de Sitter space} $\SS^{2,1}$ or \emph{anti-de Sitter space} $\H^{2,1}$, because the stabilizer of the cone is always $\SO_0(2,1)$. 
In the case of $\H^{2,1}$, we get a similar map $\dT_{g,n}\to\M^{-1}_{g,n}$, while some caution is needed in the case of $\SS^{2,1}$.

A way to understand this is to keep the same $F$ and $\EP(\ve{h},c)$ as above, but put them in $\H^{2,1}$ or $\SS^{2,1}$ as follows (see Figure \ref{figure_Sitter}): 
\begin{figure}[h]
	\centering
	\includegraphics[width=12cm]{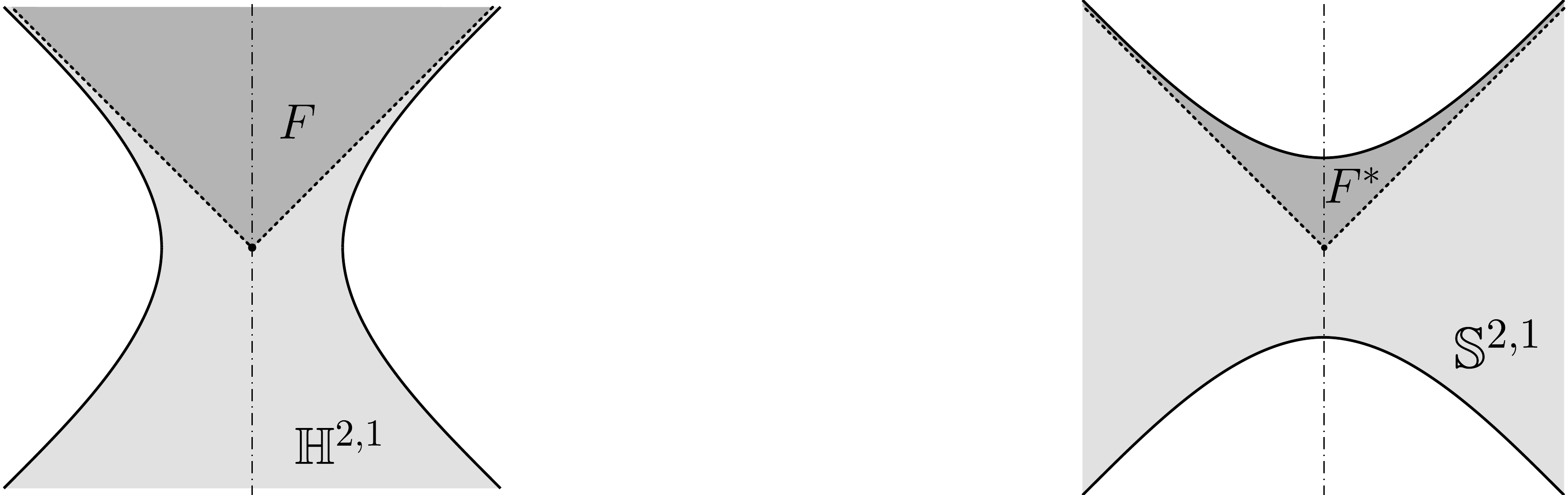}
	\caption{Section view of the future cone of the origin in the anti-de Sitter space $\H^{2,1}$ (left) and de Sitter space $\SS^{2,1}$ (right). Revolution about the vertical axis gives the actual 3d picture.}
	\label{figure_Sitter}
\end{figure}
We identify the region $\{x\pro x<1\}$ (resp.\ $\{x\pro x>-1\}$) in $\R^{2,1}$ with the projective model of $\H^{2,1}$ (resp.\ $\SS^{2,1}$) minus a plane, in which the above $F$ (resp.\ the subset $F^*:=F\cap\{x\pro x> -1\}$) is the future cone of the origin $0$. Then the surface $\pa^+\!\EP(\ve{h},c)$ can always be viewed as embedded in $\H^{2,1}$, whereas it can also be viewed as embedded in $\SS^{2,1}$ if it is contained in $F^*\subset F$, or in other words, if $(\ve{h},c)$ belongs to
$$
\dT_{g,n}^*:=\big\{(\ve{h},c)\in\dT_{g,n}\,\big|\,\text{the interior of $\EP(\ve{h},c)$ contains the hyperboloid $\H^2$}\big\}.
$$
In both cases, the polygonal pieces of $\pa^+\!\EP(\ve{h},c)$ are still spacelike, so the Lorentzian metric of $\H^{2,1}$ or $\SS^{2,1}$ induces a cone-metric $\ve{d}$ on $S_{g,n}$ of curvature $-1$ or $1$ analogous to the Euclidean one from the case of $\R^{2,1}$. We call $\ve{d}$ the hyperbolic/spherical/Euclidean \emph{Epstein-Penner metric} of $(\ve{h},c)$ in the three cases respectively.

Any spherical Epstein-Penner metric has the nontrivial property that it admits a triangulation into \emph{convex} spherical triangles (see \S \ref{subsec_21}; the similar property for $K=-1$ or $0$ is trivial, as it is satisfied by all $\ve{d}\in\M^{K}_{g,n}$).
Letting $\M^{1,*}_{g,n}$ denote the space of all $\ve{d}\in\M^1_{g,n}$ with this property, our main result is that the Epstein-Penner metric assigning maps from $\dT_{g,n}$ to $\M^{0}_{g,n}$ or $\M^{-1}_{g,n}$ and from $\dT_{g,n}^*$ to $\M^{1,*}_{g,n}$ are bijective and send the fibers of the projection $\Pi$ to \emph{discrete conformal classes}:
\begin{thmx}\label{thm_intro}
Given $g\geq0$ and $n\geq1$ with $2-2g-n<0$, every $\ve{d}\in\M^K_{g,n}$ with $K=-1$ or $0$ is the hyperbolic or Euclidean Epstein-Penner metric of a unique $(\ve{h},c)\in\dT_{g,n}$, whereas every $\ve{d}\in\M^{1,*}_{g,n}$ is the spherical Epstein-Penner metric of a unique $(\ve{h},c)\in\dT^*_{g,n}$. Moreover, distinct metrics $\ve{d},\ve{d}'$ in $\M^K_{g,n}$ or $\M^{1,*}_{g,n}$ are discretely conformal if and only if they correspond to the same $\ve{h}\in\TT_{g,n}$ with different decorations.  
\end{thmx}
The notion of discrete conformality here is introduced in  \cite{Gu-Luo-Sun-Wu_I,Gu-Guo-Luo-Sun-Wu_II} as equivalence relations in $\M^{-1}_{g,n}$ and $\M^0_{g,n}$ analogous to the usual conformality of Riemannian metrics. It naturally extends to $\M^{1,*}_{g,n}$, as recently studied in \cite{Izmestiev-Prosanov-Wu}. This is also the framework that we adopt in the present paper for spherical metrics, even through insisting working on $\M^{1,*}_{g,n}$ has some drawbacks (e.g.\ convexity of spherical triangles is not preserved by M\"obius transformations) and one can alternatively discuss discrete uniformization of spheres via Euclidean metrics as in \cite{Springborn}.

Theorem \ref{thm_intro} can be viewed as a new $3$-dimensional interpretation of discrete conformality in addition to the one of Bobenko, Pinkall and Springborn \cite{Bobenko-Pinkall-Springborn}. The latter says that $\ve{d}_1,\ve{d}_2\in\M^{-1}_{g,n},\M^0_{g,n}$ or $\M^{1,*}_{g,n}$ are discretely conformal if and only if certain $3$-dimensional polyhedral hyperbolic cone-manifolds $P(\ve{d}_1)$ and $P(\ve{d}_2)$ associated to them have the same boundary metric (see \S \ref{subsec_relationship polyhedral}), or in other words, $\pa P(\ve{d}_1)$ and $\pa P(\ve{d}_2)$ are the same hyperbolic surface with different decorations (see \S \ref{subsec_relationship decorated}). In fact, our interpretation and this one give homeomorphisms essentially inverse to each other (see Proposition \ref{prop_inverse}):
\begin{equation}\label{eqn_introinverse}
\M^K_{g,n}\xrightleftarrows[\text{Epstein-Penner metric}]{\text{Bobenko-Pinkall-Spingborn}}\dT_{g,n}
\ \text{ for $K=-1$ or $0\,,$}\ \text{ and similarly }\ 
\M^{1,*}_{g,n}\rightleftarrows\dT^*_{g,n}~,
\end{equation}
such that in each of the three cases, discrete conformal classes on the left-hand side correspond to fibers of the projection $\Pi$ on the right-hand side. This makes the study of discrete conformality equivalent to that of the decorated Teichm\"uller space $\dT_{g,n}$. An advantage of our Lorentzian construction is that it pinpoints the exact part of $\dT_{g,n}$ related to discrete conformality of spherical metrics, namely the subspace $\dT^*_{g,n}$. 
In particular, the obvious connectedness property of the fibers of $\dT^*_{g,n}$ (see Remark \ref{remark_fiberconnected}) clarifies a technical point raised in \cite[\S 2.5]{Izmestiev-Prosanov-Wu}.
However, an disadvantage of the construction is that the energy function in the variational principle (see Remark \ref{remark_HE}) does not seem to have a clear geometric meaning from the Lorentzian perspective. Also note that the Euclidean case of the construction is known to experts (see e.g.\ \cite{Bobenko-Lutz,Lutz}), so our contribution mainly lies in the hyperbolic and spherical cases.

In the development of discrete conformality, it came as a surprise that the \emph{discrete uniformization theorems} \cite{Gu-Luo-Sun-Wu_I, Gu-Guo-Luo-Sun-Wu_II} (see Theorem \ref{thm_DUT}) and the polyhedral realization theorems of Schlenker and Fillastre  \cite{Fillastre_fourier,Fillastre_dedicata, Schlenker_metriques, Schlenker_unpublished} turn out to be equivalent in certain cases via the Bobenko-Pinkall-Springborn interpretation (see Remark \ref{remark_Fillastre}).
Now, via our new interpretation, the discrete uniformization theorems become equivalent to the following results on Epstein-Penner metrics, which were not previously known:
\begin{corx}\label{coro_intro}
Given $g\geq0$ and $n\geq1$ with $2-2g-n<0$, pick $\kappa=(\kappa_i)\in(-\infty,2\pi)^n$ and let $\ve{h}\in\TT_{g,n}$ be a finite-area complete hyperbolic metric on $S^\circ_{g,n}$.
\begin{enumerate}[label=(\arabic*)]
        \item\label{item_coro1} If $\frac{1}{2\pi}\sum_{i=1}^n\kappa_i>2-2g$, then $\ve{h}$ admits a unique decoration $c$ such that the hyperbolic Epstein-Penner metric $\ve{d}\in\M^{-1}_{g,n}$ of $(\ve{h},c)$ has singular curvature $\kappa$.
        \item\label{item_coro2} If $\frac{1}{2\pi}\sum_{i=1}^n\kappa_i=2-2g$,  then $\ve{h}$ admits a decoration $c$ such that the Euclidean Epstein-Penner metric $\ve{d}\in\M^{0}_{g,n}$ of $(\ve{h},c)$ has singular curvature $\kappa$. Moreover, $c=(c_1,\cdots,c_n)$ is unique up to shrinking every $c_i$ by the same signed distance\footnote{In this paper,  ``shrinking a horocycle $c_0$ by a signed distance of $s\in\R$'' means replacing $c_0$ by the concentric horocycle in distance $|s|$ from $c_0$ which is surrounded by $c_0$ if $s>0$ and surrounds $c_0$ if $s<0$. Under the duality in Figure \ref{figure_duality}, this corresponds to scaling a point in $\pa F\setminus\{0\}$ by a factor of $e^s$.}.
        \item\label{item_coro3}  If $g=0$ (hence $n\geq3$) and either
        \begin{enumerate}[label=(\roman*)]			
        	\item\label{item_coro31} $k=(0,\cdots,0)$ and $(S^\circ_{0,n},\ve{h})$ is not the double of an ideal $n$-gon in $\H^2$, or
        	\item\label{item_coro32} $\kappa\in(0,2\pi)^n$, $\frac{1}{2\pi}\sum_{i=1}^n\kappa_i<2$ and
        	$\kappa_i<\sum_{j\neq i}\kappa_j$ for all $i$,
        \end{enumerate}
        then $\ve{h}$ admits a decoration $c$ such that $(\ve{h},c)$ belongs to $\dT^*_{g,n}$ and its spherical Epstein-Penner metric $\ve{d}\in\M^{1,*}_{0,n}$ has singular curvature $\kappa$. Moreover, $c$ is unique in Case \ref{item_coro32}.
\end{enumerate}
\end{corx}
Here, the \emph{singular curvature} $(\kappa_i)\in(-\infty,2\pi)^n$ of a metric $\ve{d}\in\M^K_{g,n}$ is defined as $\kappa_i:=2\pi-\theta_i$, where $\theta_i$ is the cone angle of $\ve{d}$ at the $i$th marked point.
It satisfies the Gauss-Bonnet formula
$$
\sum_{i=1}^n\kappa_i+K\cdot\mathrm{Area}(\ve{d})=2\pi(2-2g).
$$
This implies that the (in)equalities in Parts \ref{item_coro1} and \ref{item_coro2} are not only sufficient, but also necessary for the desired decorations to exist.  Assuming $g=0$ and $\kappa\in(0,2\pi)^n$, the two inequalities in Case \ref{item_coro32} of Part \ref{item_coro3} are necessary as well by a theorem of Luo and Tian \cite{Luo-Tian} (see also \cite{deBorbon-Panov}) on prescription of singular curvature in a usual conformal class in $\M^1_{0,n}$. In fact, the recent result of Izmestiev, Prosanov and Wu \cite{Izmestiev-Prosanov-Wu} which implies Case \ref{item_coro32} of \ref{item_coro3} is a discrete analogue of the Luo-Tian theorem.
Yet Part \ref{item_coro3} is quite limited and leaves some loose ends for future research: it would be interesting to extend the result to more general singular curvature $\kappa\in[0,2\pi)^n$ or $\kappa\in(-\infty,2\pi)^n$, or to genus $g\geq1$, even though in general we cannot expect uniqueness (see Remark \ref{remark_loosend}).

\begin{figure}[h]
	\centering
	\includegraphics[width=8.1cm]{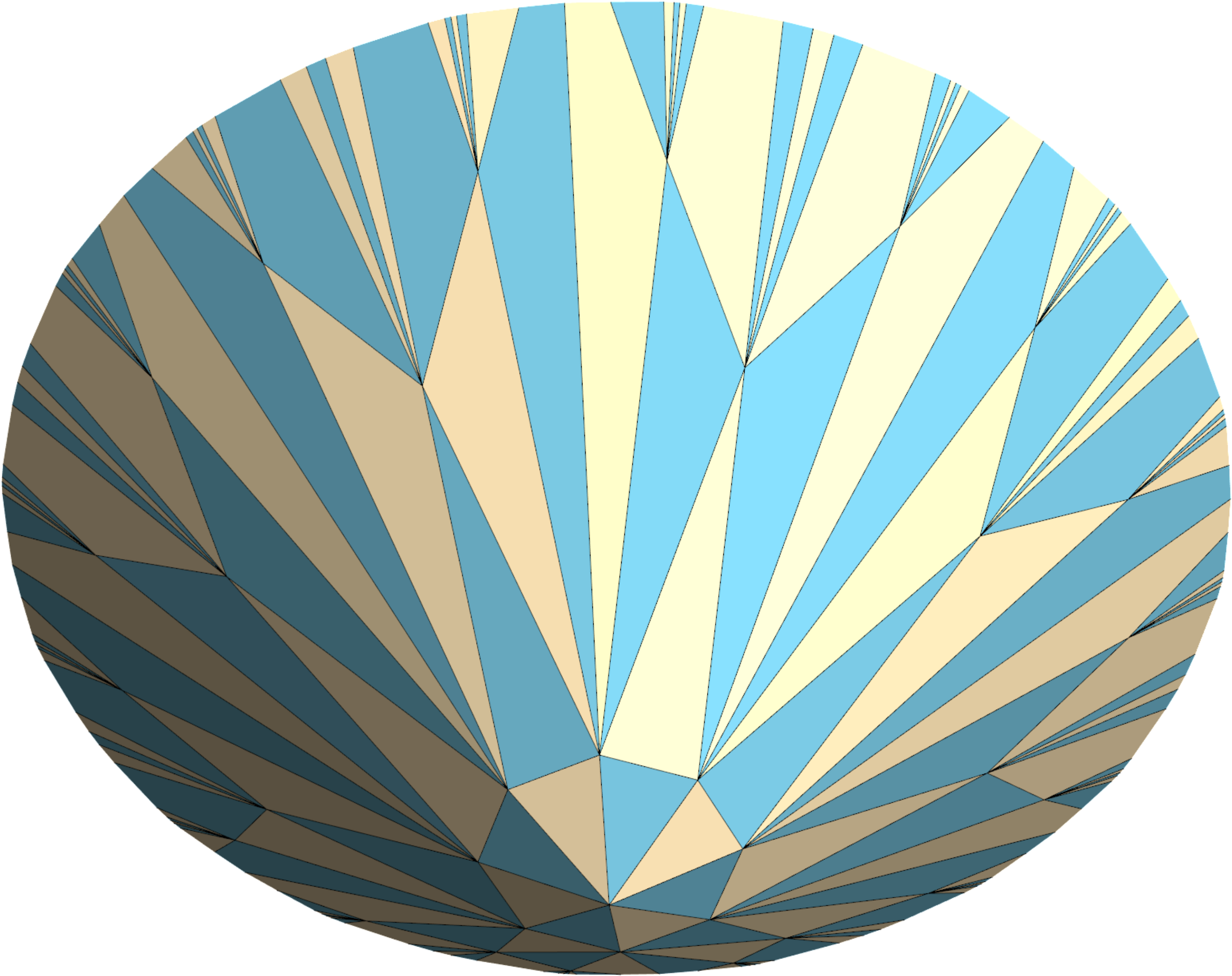}
\hspace{1.8cm}
\includegraphics[width=6.1cm]{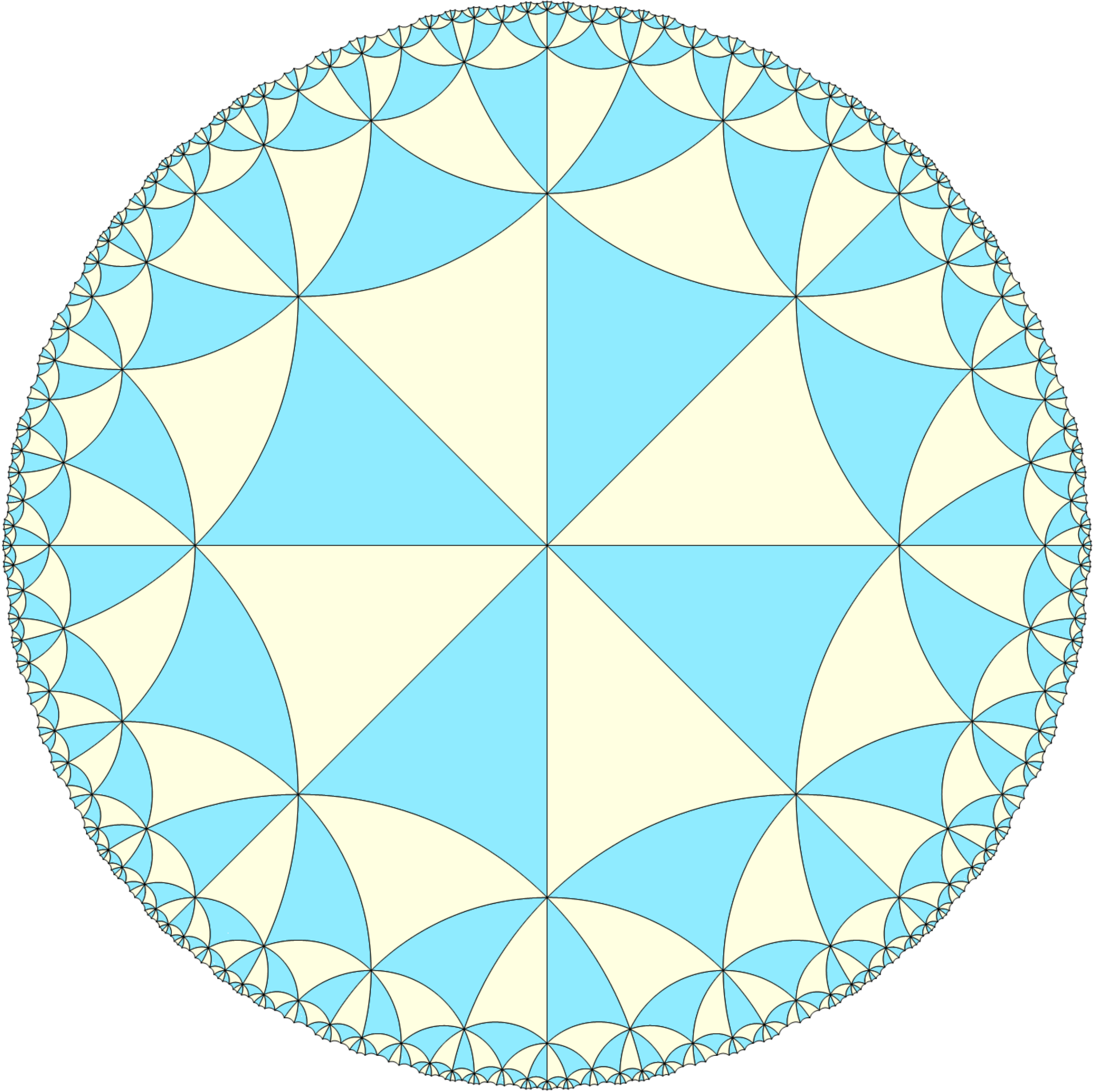}
	\caption{Boundary of a Lorentzian Fuchisan convex polyhedron and its image in $\H^2$.}
	\label{figure_FL}
\end{figure}

Finally, Corollary \ref{coro_intro} can be compared to the work of Fillastre \cite{Fillastre_annalen} on realization of cone-surfaces by \emph{Lorentzian Fuchsian convex polyhedra} in $\R^{2,1}$, $\H^{2,1}$ or $\SS^{2,1}$. To get such a polyhedron, we pick finitely many orbits in the interior of $F$ by the action of a cocompact lattice $\Gamma<\SO_0(2,1)$ (as opposed to orbits on $\pa F$ by a non-cocompact lattice in the Epstein-Penner construction) and then take their convex hull. See Figure \ref{figure_FL}.

\subsection*{Organization of the paper} We first give in \S \ref{sec_2} backgrounds on discrete conformality and the Bobenko-Pinkall-Springborn interpretation; then in \S \ref{sec_3}, after some preliminary discussions on $\H^{2,1}$ and $\SS^{2,1}$,  we prove Theorem \ref{thm_intro} and show that the maps in \autoref{eqn_introinverse}, denoted below by $\wt{\Phi}_K$ and $\Psi_K$, are inverse to each other. Crucial facts about the Lorentzian construction which lead to these results are mostly contained in \S \ref{subsec_F}.

\subsection*{Acknowledgments} The author is grateful to Tianqi Wu for enlightening comments and discussions, and to Xu Xu for his interest and the references he pointed out.

\section{Discrete conformality}\label{sec_2}
\subsection{Delaunay decomposition and triangulation}\label{subsec_21}
Recall from Introduction that $\M_{g,n}^K$ denotes the space of equivalence classes of metrics on $S_{g,n}$ of constant curvature $K\in\{-1,0,1\}$ which are allowed to have conical singularities at the $n$ marked points. 
By a \emph{triangulation} (resp.\ \emph{polygonal decomposition}) $\T$ of a metric $\ve{d}$ in $\M_{g,n}^K$, we mean a decomposition of $S_{g,n}$ into triangles\footnote{In this paper,  \emph{triangle} or \emph{polygon}, both usually denoted by $\Delta$, refers to a contractible compact subset of $\H^2$, $\E^2$ or $\SS^2$ whose boundary consists of three or at least three geodesic segments. When considering the analogous but different notion of \emph{ideal} triangle or polygon in $\H^2$, we say so explicitly and usually denote them by $\Delta'$. 
} (resp.\ polygons) given by a set $\Edge_\T$ of \emph{edges}, which are geodesic segments under the metric $\ve{d}$ joining marked points. Also denote the set of faces, namely the triangles (resp.\ polygons), by $\Face_\T$. 

We mainly consider polygons in $\H^2$, $\E^2$ or $\SS^2$ which are \emph{convex} and \emph{cyclic} in the following sense:
\begin{itemize}
\item A polygon $\Delta$ is \emph{cyclic} if it admits a circumdisk (namely a round disk $D$ in $\H^2$, $\E^2$ or $\SS^2$ such that the vertices of $\Delta$ are on the circle $\pa D$). 
\item A polygon $\Delta$ in $\SS^2$ is \emph{convex} if any two points of $\Delta$ are not antipodal and the shorter arc of $\SS^2$ joining them is in $\Delta$.
When $\Delta$ is cyclic, this is equivalent to the condition that its circumdisk has radius less than $\frac{\pi}{2}$. For polygons in $\H^2$ or $\E^2$, convexity is defined in the usual way and is implied by cyclicity.  
\end{itemize}
In particular, all triangles in $\E^2$ and $\SS^2$ are cyclic. In contrast, for a triangle in $\H^2$, there exists a unique curve of constant geodesic curvature (i.e.\ circle, horocycle or  equidistance curve to a geodesic) passing through the three vertices, and the triangle is cyclic exactly when this curve is a circle (cf.\ Remark \ref{remark_loose} below).

\begin{figure}[h]
	\centering
	\includegraphics[width=15cm]{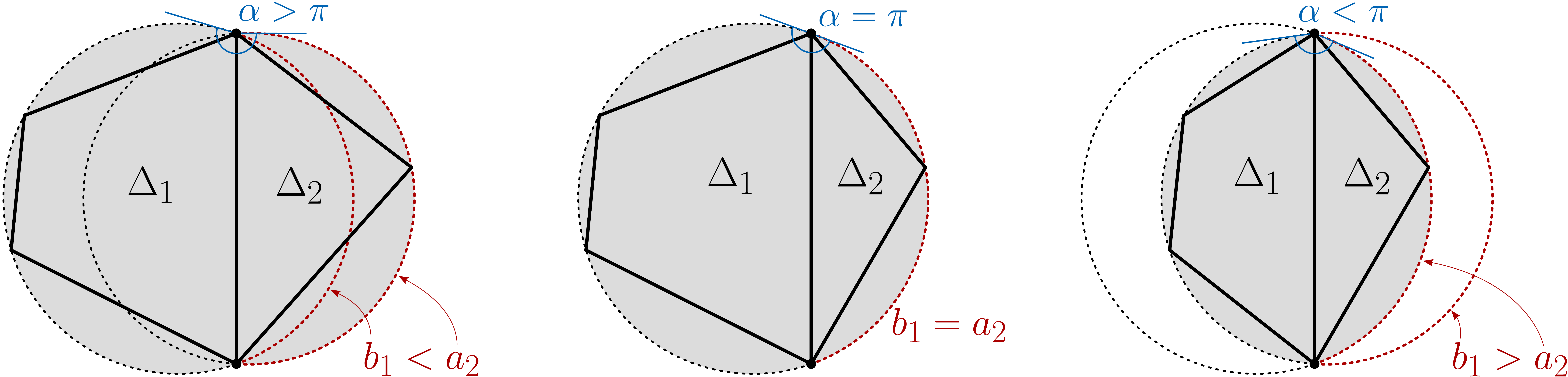}
	\caption{Adjacent convex cyclic polygons satisfying the strict Delaunay condition (left), satsifying the Delaunay condition but not the strict one (middle) and violating the Delaunay condition (right). 
The (in)equalities give equivalent characterizations of the three cases (see the proofs of Lemmas \ref{lemma_Delaunay} and \ref{lemma_FDelaunay}).
	}
	\label{figure_Delaunay}
\end{figure}

\begin{definition}\label{def_Delaunay}
A \emph{Delaunay triangulation} of a metric $\ve{d}\in \M_{g,n}^K$ is a triangulation $\T$ such that
\begin{itemize}
\item 
every face $\Delta\in\Face_\T$, developed as a triangle in $\H^2$, $\E^2$ or $\SS^2$, is convex and cyclic;
\item 
for every edge $e\in\Edge_\T$, the faces on the two sides of $e$,
  developed as adjacent triangles $\Delta_1$ and $\Delta_2$ in $\H^2$, $\E^2$ or $\SS^2$, satisfy the \emph{Delaunay condition}: the circumdisk of $\Delta_2$ does not contain any vertex of $\Delta_1$ in its interior (see Figure \ref{figure_Delaunay}).
\end{itemize}
On the other hand, a \emph{Delaunay decomposition} of $\ve{d}$ is a polygonal decomposition $\T$ such that 
\begin{itemize}
	\item
every face $\Delta\in\Face_\T$, developed as a polygon in $\H^2$, $\E^2$ or $\SS^2$, is convex and cyclic;
\item for every edge $e\in\Edge_\T$, the faces on the two sides of $e$,
developed as adjacent polygons $\Delta_1$ and $\Delta_2$ in $\H^2$, $\E^2$ or $\SS^2$, satisfy the \emph{strict Delaunay condition}: the circumdisk of $\Delta_1$ does not contain any vertex of $\Delta_2$ except the endpoints of the common edge $\Delta_1\cap\Delta_2$.
\end{itemize}
\end{definition}
Let  $\T$ be a Delaunay triangulation. Note that $\T$ is not necessarily a Delaunay decomposition, and it is if the adjacent triangles $\Delta_1,\Delta_2$ at any $e\in\Edge_\T$ satisfy the strict Delaunay condition (however, see Remark \ref{remark_generic}). We will simply call such $e$ a \emph{strictly Delaunay} edge, and otherwise call $e$ a \emph{non-strictly Delaunay} edge of $\T$. In the latter case, $\Delta_1$ and $\Delta_2$ form a convex cyclic quadrilateral in $\H^2$, $\E^2$ or $\SS^2$ split along a diagonal. We can \emph{flip} this diagonal to get a new Delaunay triangulation wherein the new edge is still non-strictly Delaunay.

The main properties of Delaunay triangulations/decompositions are summarized in the following theorem. Recall that $\M^{1,*}_{g,n}$ is the set of $\ve{d}\in\M^1_{g,n}$ admitting triangulations into convex spherical triangles.
\begin{theorem}\label{thm_Delaunay}
Every metric $\ve{d}$ in $\M^{-1}_{g,n},\M^0_{g,n}$ or $\M^{1,*}_{g,n}$ has a unique Delaunay decomposition, whereas the Delaunay triangulations of $\ve{d}$ are exactly the triangulations obtained from this decomposition by subdividing every non-triangular face into triangles. As a consequence, erasing all non-strictly Delaunay edges of any Delaunay triangulation yields the Delaunay decomposition, and moreover we can transform one Delaunay triangulation into any other by a sequence of flips along non-strictly Delaunay edges.
\end{theorem}
The Euclidean and hyperbolic cases are well known adaptations of the foundational work of Delaunay \cite{Delaunay}, see e.g.\ \cite{Gu-Luo-Sun-Wu_I,Gu-Guo-Luo-Sun-Wu_II,MR1846934,Lutz}. 
The spherical case is treated in \cite{Izmestiev-Prosanov-Wu}, where the main new ingredient is the \emph{existence} of Delaunay triangulation for $\ve{d}\in\M^{1,*}_{g,n}$.

\begin{remark}\label{remark_generic}
For a generic $\ve{d}$ in  $\M^{-1}_{g,n},\M^0_{g,n}$ or $\M^{1,*}_{g,n}$, the Delaunay decomposition of $\ve{d}$ is actually a triangulation, hence is the unique Delaunay triangulation of $\ve{d}$. This is essentially equivalent to a similar fact about $\dT_{g,n}$ concerning the Penner decomposition (see \S \ref{subsec_relationship decorated}).
\end{remark}

\begin{remark}\label{remark_loose}
We could have loosened the definition of cyclic polygons in $\H^2$, and hence the Delaunay condition, by allowing \emph{generalized disks} (i.e.\ convex domains in $\H^2$ whose boundaries are curves of constant geodesic curvature) as circumdisks. This would make every triangle cyclic. However, it does not change anything because by \cite[Thm.\ 14]{Gu-Guo-Luo-Sun-Wu_II} (see also \cite[Lemma 3.6]{Prosanov}), for a metric in $\M^{-1}_{g,n}$, any face of a Delaunay decomposition or triangulation under this looser definition is actually cyclic in our sense.
\end{remark}

\subsection{Discrete conformality and uniformization}\label{subsec_DUT}
The following definition is introduced in \cite{Gu-Luo-Sun-Wu_I, Gu-Guo-Luo-Sun-Wu_II} for $\M^{0}_{g,n}$ and $\M^{-1}_{g,n}$. The natural extension to $\M^{1,*}_{g,n}$ is studied in  \cite{Izmestiev-Prosanov-Wu}.
\begin{definition}\label{def_DC}
Two metrics $\ve{d}$ and $\ve{d}'$ on $S_{g,n}$ both in $\M_{g,n}^{-1}$, $\M_{g,n}^0$ or $\M_{g,n}^{1,*}$ are said to be \emph{discretely conformal} if there exists a sequence $\ve{d}=\ve{d}_0,\ve{d}_1,\ve{d}_2,\cdots,\ve{d}_N=\ve{d'}$ in the same space and a Delaunay triangulation $\T_r$ of each $\ve{d}_r$ such that every $(\ve{d}_{r+1},\T_{r+1})$ is related to $(\ve{d}_r,\T_r)$ in either of the following two ways:
\begin{itemize}
	\item \emph{flips}: $\ve{d}_{r+1}=\ve{d}_r$, but $\T_{r+1}$ and $\T_r$ are different Delaunay triangulations of this metric. By Theorem \ref{thm_Delaunay}, this means that $\T_{r+1}$ is obtained from $\T_r$ by a sequence of flips along non-strictly Delaunay edges.
	\item \emph{vertex scaling}: $\T_{r+1}$ and $\T_r$ are isotopic and there exists $(u_i)\in\R^n$ such that
    for any $i,j\in\{1,\cdots,n\}$ and any edge $e$ joining the $i$th and $j$th marked points, the lengths of $e$ under the metrics $\ve{d}_r$ and $\ve{d}_{r+1}$, denoted respectively by $\ell_e$ and $\hat{\ell}_e$, satisfy the relation
	\begin{equation}\label{eqn_vertexscaling}
	\begin{cases}
		\sinh\tfrac{\hat\ell_e}{2}=e^{\frac{u_i+u_j}{2}}\sinh\tfrac{\ell_e}{2}&\text{ for $\M_{g,n}^{-1}$,}\\[5pt]    
		\hat\ell_e=e^{\frac{u_i+u_j}{2}}\ell_e&\text{ for $\M_{g,n}^0$,}\\[5pt]
		\sin\tfrac{\hat\ell_e}{2}=e^{\frac{u_i+u_j}{2}}\sin\tfrac{\ell_e}{2}&\text{ for $\M_{g,n}^{1,*}$.}
	\end{cases}
	\end{equation}
\end{itemize}
This defines an equivalence relation on each of the spaces $\M_{g,n}^{-1}$, $\M_{g,n}^0$ and $\M_{g,n}^{1,*}$. The equivalence classes are call \emph{discrete conformal classes}. 
\end{definition}
\begin{example}\label{example_double}
For $g=0$ and $n\geq3$, a particular type of metric $\ve{d}$ in $\M^{-1}_{0,n}$, $\M^0_{0,n}$ or $\M^{1,*}_{0,n}$ is \emph{doubled $n$-gon}. We call the cone-surface $(S_{0,n},\ve{d})$ a doubled $n$-gon if it can be identified homeomorphically with the sphere such that the $n$ marked points are on the equator and the reflection about the equator, which switches the north and south hemispheres, is an isometry. In other words, this means that $(S_{0,n},\ve{d})$ is the double of a (not necessarily convex) polygon in $\H^2$, $\E^2$ or $\SS^2$. Abusing the terminology, we also call the metric $\ve{d}$ a doubled $n$-gon in this case.
Given such a $\ve{d}$, it can be shown that any metric discretely conformal to $\ve{d}$ is again a doubled $n$-gon, and there exists a metric discretely conformal to $\ve{d}$ which is the double of a convex cyclic $n$-gon\footnote{This follows e.g.\ from the interpretation of discrete conformality through $\dT_{g,n}$ explained later (namely the bijection \autoref{eqn_introinverse}), basically because $\ve{d}$ is a doubled $n$-gon if and only if the corresponding $(\ve{h},c)\in\dT_{g,n}$ is the double of a decorated ideal $n$-gon in $\H^2$, while we can always modify a decorated ideal $n$-gon into a \emph{cyclic} one (in the sense specified in \S \ref{subsec_relationship decorated} below) by a change of decoration.}.
\end{example}
In both $\M^0_{g,n}$ and $\M^{1,*}_{g,n}$, there are natural equivalence relations finer than discrete conformality. The one in $\M^0_{g,n}$ just comes from \emph{dilation}, namely multiplying a Euclidean metric by a positive constant. The one in $\M^{1,*}_{g,n}$ is more subtle: it is restricted from the \emph{M\"obius equivalence} in $\M^1_{g,n}$, where two spherical metrics $\ve{d}_1,\ve{d}_2\in\M^1_{g,n}$ are said to be M\"obius-equivalent if $\ve{d}_i$ has a developing pair $(f_i,\rho_i)$ (where $\rho_i:\pi_1(S^\circ_{g,n})\rightarrow\SO(3)$ is a \emph{holonomy representation}, and the \emph{developing map} $f_i$ is a $\rho_i$-equivariant local diffeomorphism from the universal cover of $S^\circ_{g,n}$ to $\SS^2$; see \cite{Goldman}) such that $(f_1,\rho_1)$ is \emph{conjugate} to $(f_2,\rho_2)$ through a  M\"obius transformation $a\in\PSL(2,\C)$, in the sense that $f_2=a\circ f_1$ and $\rho_2(\gamma)=a\circ \rho_1(\gamma)\circ a^{-1}$ for all $\gamma\in\pi_1(S^\circ_{g,n})$\footnote{If $a\in\SO(3)\subset\PSL(2,\C)$, then this conjugation just means that $\ve{d}_1$ and $\ve{d}_2$ are the same point in $\M^1_{g,n}$. A ``nontrivial'' conjugation, through some $a\notin\SO(3)$, occurs only when $\rho_1$ and $\rho_2$ are either both trivial or both \emph{co-axial} (a representation in $\SO(3)$ is said to be co-axial if its image is in a one-parameter subgroup of rotations). In other words, the M\"obius equivalence class of $\ve{d}\in\M^1_{g,n}$ consists not solely of $\ve{d}$ itself only when the holonomy of $\ve{d}$ is either trivial or co-axial, see \cite[\S 2]{Eremenko_survey}.}.

The vertex scaling operation first appeared without the Delaunay condition and was originally taken as definition for discrete conformality (see \cite{Bobenko-Pinkall-Springborn,Luo_combinatorial}). It is a natural analogue of the usual conformality for Riemannian metrics. 
The legitimacy of the above more sophisticated definition is justified by the following:

\begin{theorem}[\textbf{Discrete Uniformization}]\label{thm_DUT}
Given $g\geq0$ and $n\geq1$ with $2-2g-n<0$, pick $\kappa=(\kappa_i)\in(-\infty,2\pi)^n$.
	\begin{enumerate}[label=(\arabic*)]
		\item\label{item_DUT1} If $\frac{1}{2\pi}\sum_{i=1}^n\kappa_i>2-2g$, then every discrete conformal class in $\M^{-1}_{g,n}$ contains a unique element with singular curvature $\kappa$.
		\item\label{item_DUT2} If $\frac{1}{2\pi}\sum_{i=1}^n\kappa_i=2-2g$, then every discrete conformal class in $\M^{0}_{g,n}$ contains an element with singular curvature $\kappa$, and this element is unique up to dilation.
	\end{enumerate}	
	\vspace{-0.15cm}
	\begin{enumerate}[label=(3.\roman*)]
		\item\label{item_DUT3} If $g=0$ and $\kappa=(0,\cdots,0)$, then every discrete conformal class in $\M_{0,n}^{1,*}$ not formed by doubled $n$-gons (see Example \ref{example_double}) contains an element with singular curvature $\kappa$, and this element is unique up to M\"obius equivalence.
		\item\label{item_DUT4}
       	If $g=0$, $\kappa\in(0,2\pi)^n$, $\frac{1}{2\pi}\sum_{i=1}^n\kappa_i<2$ and
			$\kappa_i<\sum_{j\neq i}\kappa_j$ for all $i$,
		    then every discrete conformal class in $\M^{1,*}_{0,n}$ contains a unique element with singular curvature $\kappa$. 
	\end{enumerate}
\end{theorem}
 Parts \ref{item_DUT1} and \ref{item_DUT2} are due to  Gu-Guo-Luo-Sun-Wu \cite{Gu-Guo-Luo-Sun-Wu_II} and Gu-Luo-Sun-Wu \cite{Gu-Luo-Sun-Wu_I}, respectively, while \ref{item_DUT4} is obtained recently by  Izmestiev-Prosanov-Wu \cite{Izmestiev-Prosanov-Wu}. Part \ref{item_DUT3} has not appeared explicitly in the literature, but is equivalent to Rivin's polyhedral realization theorem \cite{Rivin_intrinsic} for hyperbolic puncture spheres. In fact, the doubled $n$-gons excluded in \ref{item_DUT3} correspond to those ideal polyhedra in $\H^3$ from Rivin's theorem which degenerate to two-sided ideal polygons. Another equivalent version is given by Springborn \cite{Springborn}, where he studies discrete uniformization of spheres via Euclidean metrics instead.
 
As is common in the study of spherical metrics comparing to Euclidean and hyperbolic ones, the only discrete uniformization results in spherical setting, namely \ref{item_DUT3} and \ref{item_DUT4}, are much less complete than the Euclidean and hyperbolic counterparts, hence provides a direction for future research. See Remark \ref{remark_loosend}.
\begin{remark}
As mentioned in Introduction, by the Gauss-Bonnet formula, the (in)equalities in Parts \ref{item_DUT1} and \ref{item_DUT2} of Theorem \ref{thm_DUT} are also necessary for the desired elements in $\M^{-1}_{g,n}$ and $\M^0_{g,n}$ to exist, whereas by the Luo-Tian theorem \cite{Luo-Tian}, the two inequalities in Part \ref{item_DUT4} are necessary as well under the assumptions $g=0$ and $\kappa\in(0,2\pi)^n$. In fact, \ref{item_DUT4} is a discrete analogue of the Luo-Tian theorem.
\end{remark}

\begin{remark}\label{remark_loosend}
Parts \ref{item_DUT1} and \ref{item_DUT2} are discrete analogues of the McOwen-Troyanov uniformization theorem \cite{McOwen, Troyanov}, which asserts the unique-existence of a Euclidean or hyperbolic metric with prescribed singular curvature $\kappa$ (satisfying the Gauss-Bonnet constraint, namely the (in)equality in \ref{item_DUT1} or \ref{item_DUT2}) in a usual conformal class in $\M^0_{g,n}$ or $\M^{-1}_{g,n}$. If $\kappa=(0,\cdots,0)$, this reduces to the Poincar\'e-Koebe uniformization theorem. 
In contrast, the problem of prescribing singular curvature for spherical metrics is much more elusive: the Gauss-Bonnet constraint is not sufficient for existence, and there is no uniqueness in general. There are existence results under various assumptions (see e.g.\ \cite{Li-Song-Xu} and the references therein), among which the aforementioned Luo-Tian theorem is the most fundamental one and has a rare uniqueness part.
It would be interesting to seek discrete analgoues of other existence results, such as \cite[Thm.\ 1.1]{Bartolucci} and \cite[Thm.\ C]{Troyanov}. 
\end{remark}

\subsection{Relationship with polyhedra in $\H^3$}\label{subsec_relationship polyhedral}
Bobenko, Pinkall and Springborn \cite{Bobenko-Pinkall-Springborn} find a $3$-dimensional hyperbolic geometric interpretation of discrete conformality.
It involves a hyperbolic cone-manifold $P(\ve{d})$ with convex ideal polyhedral boundary constructed from $\ve{d}\in\M^{-1}_{g,n}$, $\M^0_{g,n}$ or $\M^{1,*}_{g,n}$. Such cone-manifolds and their analogues are the main objects of study in the variational approach to Alexandrov type polyhedral realization problems \cite{Bobenko-Izmestiev,Fillastre-Izmestiev_1,Fillastre-Izmestiev_2,Izmestiev-Prosanov-Wu,Prosanov,Springborn} and have been given various names in different settings. Here we call them \emph{generalized Fuchsian convex polyhedra}. In order to define them, we first introduce the building blocks.

\begin{figure}[h]
	\centering
	\includegraphics[width=16.6cm]{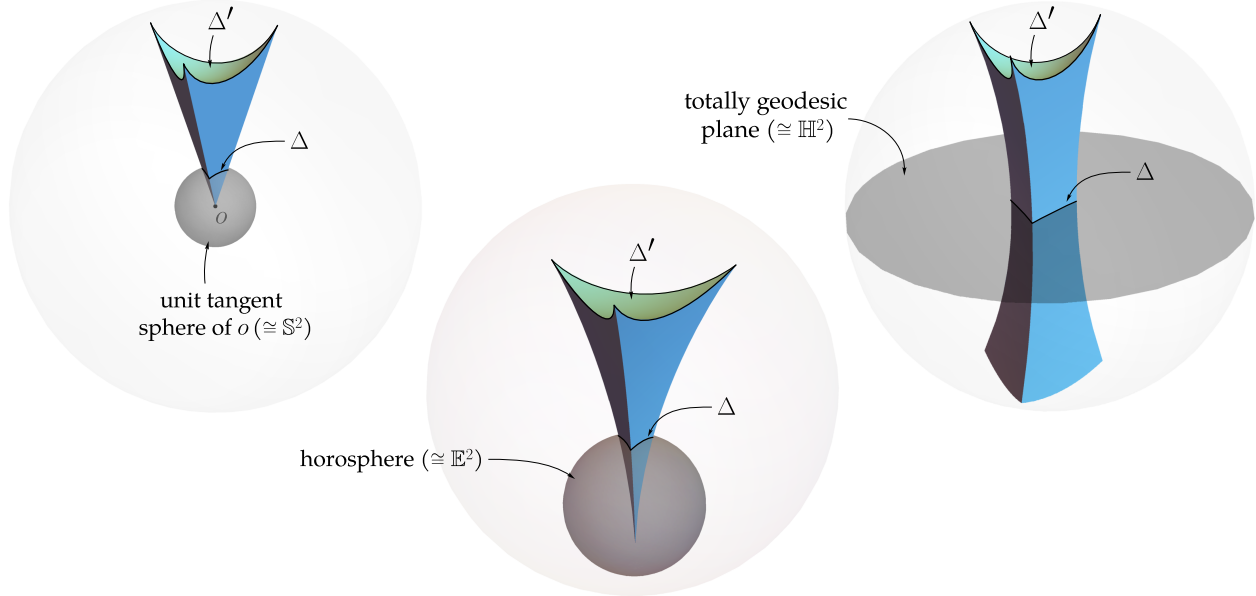}
	\caption{The polyhedron $P(\Delta)$ in the Poincar\'e ball model of $\H^3$, when $\Delta$ is a triangle in $\SS^2$, $\E^2$ and $\H^2$, respectively. The vertices of the face $\Delta'$ are ideal.}
	\label{figure_bricks}
\end{figure}

Identify $\SS^2$ with the unit tangent sphere of a point in the hyperbolic $3$-space $\H^3$, and identify $\E^2$ (resp.\ $\H^2$) with a horosphere (resp.\ a totally geodesic plane).   
Given a convex cyclic polygon $\Delta$ in $\H^2$, $\E^2$ or $\SS^2$, we let $P(\Delta)\subset\H^3$ denote the polyhedron as shown in Figure \ref{figure_bricks}, 
whose boundary is formed by the side faces passing through each edge of $\Delta$ and orthogonal to $\Delta$, along with an ideal polygonal face denoted by $\Delta'$.

\begin{figure}[h]
	\centering
	\includegraphics[width=8.5cm]{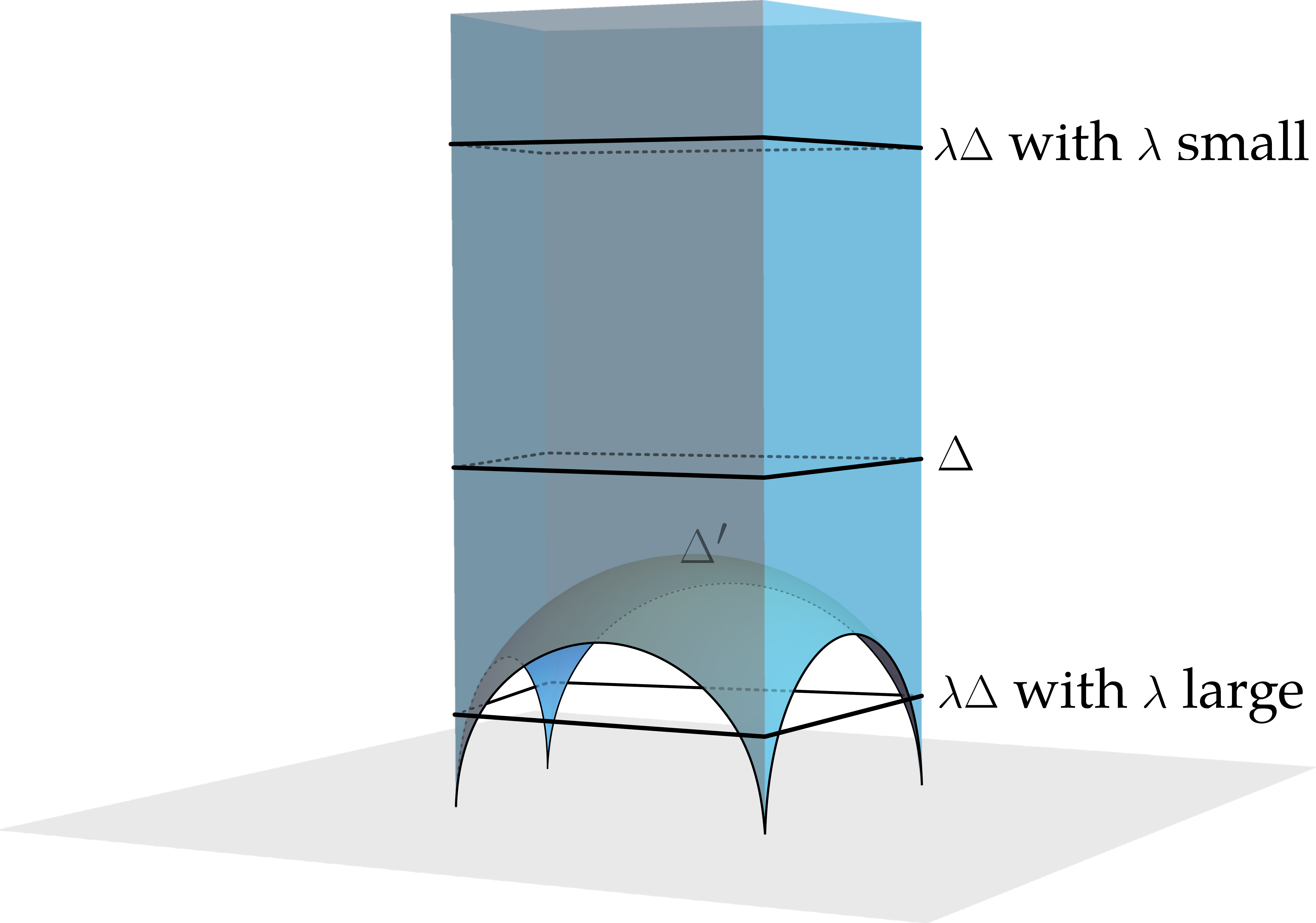}
	\caption{$P(\Delta)\cong P(\lambda\Delta)$ in the upper half-space model for a cyclic quadrilateral $\Delta\subset\E^2$.}
	\label{figure_upperhalfspace}
\end{figure}

\begin{remark}
	$P(\Delta)$ can by understood in the Klein model of $\H^3$ as a pyramid whose base is an ideal polygon, while the apex is either an ordinary point, an ideal point or a hyper-ideal point depending on whether $\Delta$ is in $\SS^2$, $\E^2$ or $\H^2$. In the case of $\H^2$, we may truncate $P(\Delta)$ along the plane polar to the apex (i.e.\ the plane containing $\Delta$) to get a polyhedron of finite volume, but this is unnecessary for our purpose.
\end{remark}
\begin{remark}\label{remark_scaling}
When $\Delta\subset\E^2$, the polyhedron $P(\Delta)$ does not necessarily contain $\Delta$. In fact, for any dilation $\lambda\Delta$ of $\Delta$ ($\lambda>0$), by varying the base horosphere in a concentric way, we may view $P(\lambda\Delta)$ and $P(\Delta)$ as the same ideal polyhedron in $\H^3$, only with the position of the base polygon $\lambda\Delta$ changing. When $\lambda$ is large, this polygon is close to the boundary of $\pa \H^3$ and intersects the face $\Delta'$. See Figure \ref{figure_upperhalfspace}.
\end{remark}
\begin{remark}
When $\Delta\subset\H^2$ is cyclic in the looser sense mentioned in Remark \ref{remark_loose} (i.e.\ its vertices are on a curve $C$ of constant geodesic curvature), $P(\Delta)$ is still well defined, and $\Delta$ is cyclic in our sense (i.e.\ $C$ is a circle) if and only if the plane containing the face $\Delta'$ is ultraparallel to the one containing $\Delta$. See \cite{Prosanov}. 
\end{remark}

\begin{figure}[h]
	\centering
	\includegraphics[width=14.1cm]{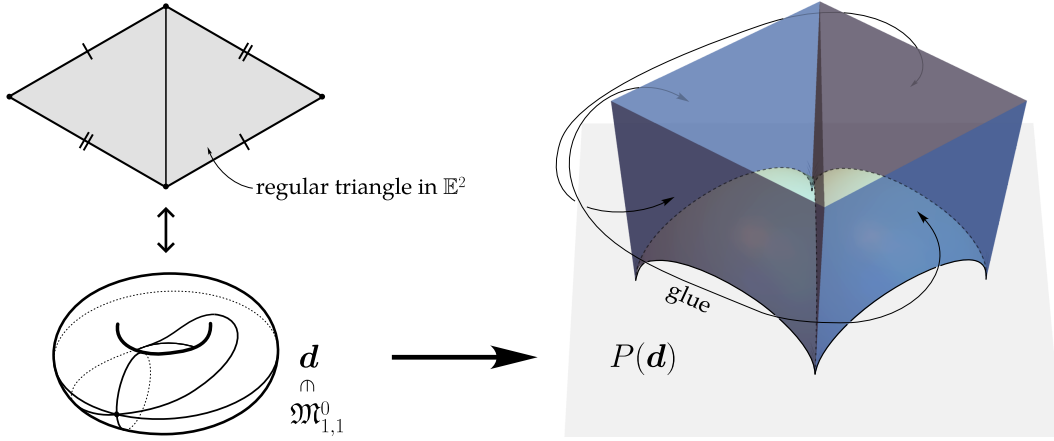}
	\caption{A simple example of $P(\ve{d})$, where $\ve{d}$ is a nonsingular Euclidean metric on the torus with one marked point and its Delaunay decomposition is a triangulation.}
	\label{figure_Pd}
\end{figure}

\begin{definition}
For any metric $\ve{d}$ in $\M^{-1}_{g,n}$, $\M^0_{g,n}$ or $\M^{1,*}_{g,n}$, the associated \emph{generalized Fuchsian convex polyhedron}, denoted by $P(\ve{d})$, is the hyperbolic $3$-manifold with boundary and with conical singular lines constructed by gluing together the polyhedra $P(\Delta)$ for all the faces $\Delta\in\Face_\T$ of the Delaunay decomposition $\T$ of $\ve{d}$ (according to the combinatorics of $\T$). Equivalently, $P(\ve{d})$ is obtained by gluing $P(\Delta)$ for all the faces $\Delta$ of any Delaunay triangulation of $\ve{d}$.  See Figure \ref{figure_Pd}. 
We call $P(\ve{d})$ a generalized Fuchsian convex polyhedron of \emph{hyperbolic}, \emph{parabolic} or \emph{elliptic} type when  $\ve{d}$ in $\M^{-1}_{g,n}$, $\M^0_{g,n}$ or $\M^{1,*}_{g,n}$, respectively.
\end{definition}	
A crucial fact behind this definition is that two adjacent convex cyclic polygons $\Delta_1$ and $\Delta_2$ satisfy the Delaunay condition (resp.\ strict Delaunay condition) if and only if the dihedral angle between the faces $\Delta'_1$ and $\Delta'_2$ of $P(\Delta_1)\cup P(\Delta_2)$ is $\leq\pi$ (resp.\ $<\pi$). See Lemma \ref{lemma_Delaunay} below. The equality case of this fact implies that the Delaunay decomposition and triangulations of $\ve{d}$ yield the same $P(\ve{d})$. Also, the fact implies that $P(\ve{d})$ does not have 
boundary dihedral angles larger than $\pi$, which explains the term ``convex''.

The boundary $\pa P(\ve{d})$ is intrinsically a finite-area hyperbolic surface constructed by assembling the face $\Delta'$ of $P(\Delta)$ for all $\Delta\in\Face_\T$. It can be shown that this surface is \emph{complete}. Meanwhile, for each $P(\Delta)$, there is a natural diffeomorphism from $\Delta$ with vertices removed to $\Delta'$, given by rays issuing orthogonally from $\Delta$. They fit together to form a diffeomorphism from $(S_{g,n}^\circ,\ve{d})$ (the cone-surface $(S_{g,n},\ve{d})$ with marked points removed) to the hyperbolic surface $\pa P(\ve{d})$. Therefore, the metric of $\pa P(\ve{d})$ pulls back to one on $S_{g,n}^\circ$ and represents a point in the Teichm\"uller space $\TT_{g,n}$. This defines maps
\begin{align*}
\Phi_{-1} \text{ (resp.\ $\Phi_0$, $\Phi_1$)}:\ \M^{-1}_{g,n} \text{ (resp.\ $\M^{0}_{g,n}$, $\M^{1,*}_{g,n}$)}&\to \TT_{g,n},\\
\ve{d}&\mapsto \left[\parbox{6.6cm}{pullback of the hyperbolic metric on $\pa P(\ve{d})$\\ by the diffeomorphism from $(S_{g,n}^\circ,\ve{d})$}\right].
\end{align*}
It can be shown that dilation-equivalent Euclidean metrics in $\M^0_{g,n}$ (resp.\ M\"obius-equivalent spherical metrics in $\M^{1,*}_{g,n}$) are sent by $\Phi_0$ (resp.\ $\Phi_1$) to the same point of $\TT_{g,n}$.

\begin{remark}\label{remark_hilbert}
The construction of the hyperbolic metric on each convex cyclic polygon $\Delta$ (and hence the construction of $\Phi_K(\ve{d})$)
can also be reformulated purely through $2$-dimensional geometry: if we identify the circumdisk $D_\Delta$ of $\Delta$ with the Klein model of $\H^2$, then the restriction of its hyperbolic metric to $\Delta$ coincides with the metric pulled back from $\Delta'$, see \cite[\S 4.1]{Bobenko-Lutz} or \cite[\S 5.1]{Bobenko-Pinkall-Springborn}. This can be shown by looking into the identification between  $D_\Delta$ and the plane in $\H^3$ containing $\Delta'$ which extends the above diffeomorphism $\Delta\overset\sim\to\Delta'$, or alternative by our Lorentzian construction in \S \ref{subsec_F}.
\end{remark}

Now the interpretation of discrete conformality in \cite{Bobenko-Pinkall-Springborn} can be stated as:
\begin{proposition}\label{prop_BPS}
The discrete conformal classes in $\M^{-1}_{g,n}$, $\M^0_{g,n}$ or $\M^{1,*}_{g,n}$ are the fibers of the map $\Phi_{-1}$, $\Phi_0$ or $\Phi_1$.
\end{proposition}
A more refined statement is given in Theorem \ref{thm_BPS} below along with the proof. The above statement allows us to reformulate Theorem \ref{thm_DUT} as the following equivalent theorem about realizability of finite-area complete hyperbolic surfaces as boundaries of generalized Fuchsian convex polyhedra:
 
 \begin{theorem}[\textbf{Discrete Uniformization, 2nd version}]\label{thm_DUT2}
Given $g\geq0$ and $n\geq1$ with $2-2g-n<0$, pick $\kappa=(\kappa_i)\in(-\infty,2\pi)^n$ and let $\ve{h}\in\TT_{g,n}$ be a finite-area complete hyperbolic metric on $S^\circ_{g,n}$.
	\begin{enumerate}[leftmargin=1cm, label=(\arabic*)]
	\item\label{item_DUT21} If $\frac{1}{2\pi}\sum_{i=1}^n\kappa_i>2-2g$, then $\ve{h}$ can be realized as the boundary metric of a unique generalized Fuchsian convex polyhedron of hyperbolic type with singular curvature $\kappa$.
	\item\label{item_DUT22} If $\frac{1}{2\pi}\sum_{i=1}^n\kappa_i=2-2g$, then $\ve{h}$ can be realized as the boundary metric of a unique generalized Fuchsian convex polyhedron of parabolic type with singular curvature $\kappa$.
	\item\label{item_DUT23} 
   If $g=0$ and either
   \begin{enumerate}[label=(\roman*)]			
   	\item\label{item_DUT231} $k=(0,\cdots,0)$ and $(S^\circ_{0,n},\ve{h})$ is not the double of an ideal $n$-gon in $\H^2$, or
   	\item\label{item_DUT232} $\kappa\in(0,2\pi)^n$, $\frac{1}{2\pi}\sum_{i=1}^n\kappa_i<2$ and
   	$\kappa_i<\sum_{j\neq i}\kappa_j$ for all $i$,
   \end{enumerate}
	then $\ve{h}$ can be realized as the boundary metric of a unique generalized  Fuchsian  convex polyhedron of elliptic type with singular curvature $\kappa$.
\end{enumerate}
\end{theorem}
Here, the singular curvature of a $3$-manifold with conical singular lines is originally defined by the cone angle at each singular line (see e.g.\ \cite{Fillastre-Izmestiev_1}). However, for the particular manifold $P(\ve{d})$, it coincides with the singular curvature $(\kappa_i)\in(-\infty,2\pi)^n$ of the base surface $(S_{g,n},\ve{d})$.
\begin{remark}\label{remark_Fillastre}
The nonsingular (i.e.\ $\kappa=(0,\cdots,0)$) case of Theorem \ref{thm_DUT2} has been known before the discovery of the relationship with discrete conformality.
In fact, Case \ref{item_DUT231} of Part \ref{item_DUT23} is exactly the theorem of Rivin mentioned after Theorem \ref{thm_DUT}, whose degenerate situation (two-sided ideal polygons in $\H^3$) corresponds to the doubled ideal $n$-gons excluded in \ref{item_DUT231}. Meanwhile, in the nonsingular case of Parts \ref{item_DUT21} and \ref{item_DUT22}, $P(\ve{d})$ is the quotient of a polyhedral convex domain in $\H^3$, with infinitely many faces, by an isometric action of $\pi_1(S_g)$ which preserves a horosphere (if $g=1$) or a plane (if $g\geq2$). Such quotient polyhedra, as well as their variants with ideal vertices replaced by non-ideal or hyper-ideal ones, have been extensively studied \cite{Fillastre_fourier,Fillastre_dedicata,Fillastre-Izmestiev_1,Schlenker_metriques, Schlenker_unpublished}. In particular, the nonsingular case of \ref{item_DUT21} and \ref{item_DUT22} are contained in Fillastre's more general result \cite{Fillastre_dedicata}.
\end{remark}

\subsection{Relationship with decorated Teichm\"uller space}\label{subsec_relationship decorated}
A closer look into the construction of the generalized Fuchsian convex polyhedron $P(\ve{d})$ reveals a link with the decorated Teichm\"uller space $\dT_{g,n}$, formulated as the following theorem, which strengthens Proposition \ref{prop_BPS}: 

\begin{theorem}\label{thm_BPS}
The map $\Phi_K$ naturally factorizes as an injective map $\wt{\Phi}_K$ from $\M^{-1}_{g,n}$, $\M^0_{g,n}$ or $\M^{1,*}_{g,n}$ to $\dT_{g,n}$ 
composed with the projection $\Pi:\dT_{g,n}\to\TT_{g,n}$. Two metrics $\ve{d}_1,\ve{d}_2$ in $\M^{-1}_{g,n}$, $\M^0_{g,n}$ or $\M^{1,*}_{g,n}$ are discretely conformal if and only if $\wt{\Phi}_K(\ve{d}_1),\wt{\Phi}_K(\ve{d}_2)\in\dT_{g,n}$ are in the same fiber of $\Pi$.
\end{theorem}	
The map $\wt{\Phi}_K$ is the left-to-right arrow in \autoref{eqn_introinverse} of Introduction.
We will show in Proposition \ref{prop_inverse} that the Epstein-Penner metric assigning map $\Psi_K$ is inverse to $\wt{\Phi}_K$, which implies that $\wt{\Phi}_K$ is bijective when $K=-1,0$ and is bijective to the subspace $\dT^*_{g,n}$ of the target $\dT_{g,n}$ when $K=1$. The bijectivity when $K=-1$ is equivalent to \cite[Cor.\ 4.2, Lemma 4.3]{Prosanov}.
\begin{remark}\label{remark_HE}
	The construction of $\wt{\Phi}_K$ appeared implicitly in \cite{Bobenko-Pinkall-Springborn,Izmestiev-Prosanov-Wu,Prosanov} and is fundamental for the variational approach to Alexandrov type polyhedral realization problems. In fact, if we identify a fiber $\Pi^{-1}(\ve{h})\subset\dT_{g,n}$ with $\R^n$ by choosing an auxiliary decoration for  $\ve{h}\in\TT_{g,n}$, then the map $\wt{\Phi}_K^{-1}(\Pi^{-1}(\ve{h}))\hookrightarrow\Pi^{-1}(\ve{h})\cong\R^n$, viewed as coordinates on $\wt{\Phi}_K^{-1}(\Pi^{-1}(\ve{h}))$, is exactly the ``distance coordinates'' for generalized Fuchsian convex polyhedra with fixed boundary metric $\ve{h}$ as considered in \cite{Izmestiev-Prosanov-Wu,Prosanov} (adapted in turn from the works \cite{Bobenko-Izmestiev,Fillastre-Izmestiev_1} on non-ideal settings).  A key object in this variational approach is the \emph{Hilbert-Einstein functional} as a function in these coordinates. However, since these coordinates are not needed in this paper, our presentation in the proof below is slightly different from the above references in that we do not pick any auxiliary decoration.
\end{remark}

Before giving the proof of Theorem \ref{thm_BPS}, we recall some backgrounds on $\dT_{g,n}$ analogous to the constructions in \S \ref{subsec_21}. An ideal polygon in $\H^2$ is said to be \emph{decorated} if a horocycle is chosen at each vertex. A decorated ideal polygon $\Delta'$ is said to be \emph{cyclic} if it satisfies either of the following equivalent conditions (cf.\ \cite[\S 4]{Springborn}): 
\begin{enumerate}[label=(\alph*)]
	\item\label{item_d1} After shrinking all the horocycles by the same distance if necessary, there exists a circle in $\H^2$ tangent to each of them.
	\item\label{item_d2} After shrinking all the horocycles by the same distance if necessary, there exists a circle in $\H^2$ intersecting each of them orthogonally.
	\item\label{item_d3} The points on the light cone $\pa F\subset\R^{2,1}$ corresponding to the horocycles (under the duality in Figure \ref{figure_duality}) lie on some affine plane $A\subset \R^{2,1}$ such that $A\cap\pa F$ is an ellipse.
\end{enumerate}
\begin{figure}[h]
	\centering
	\includegraphics[width=13cm]{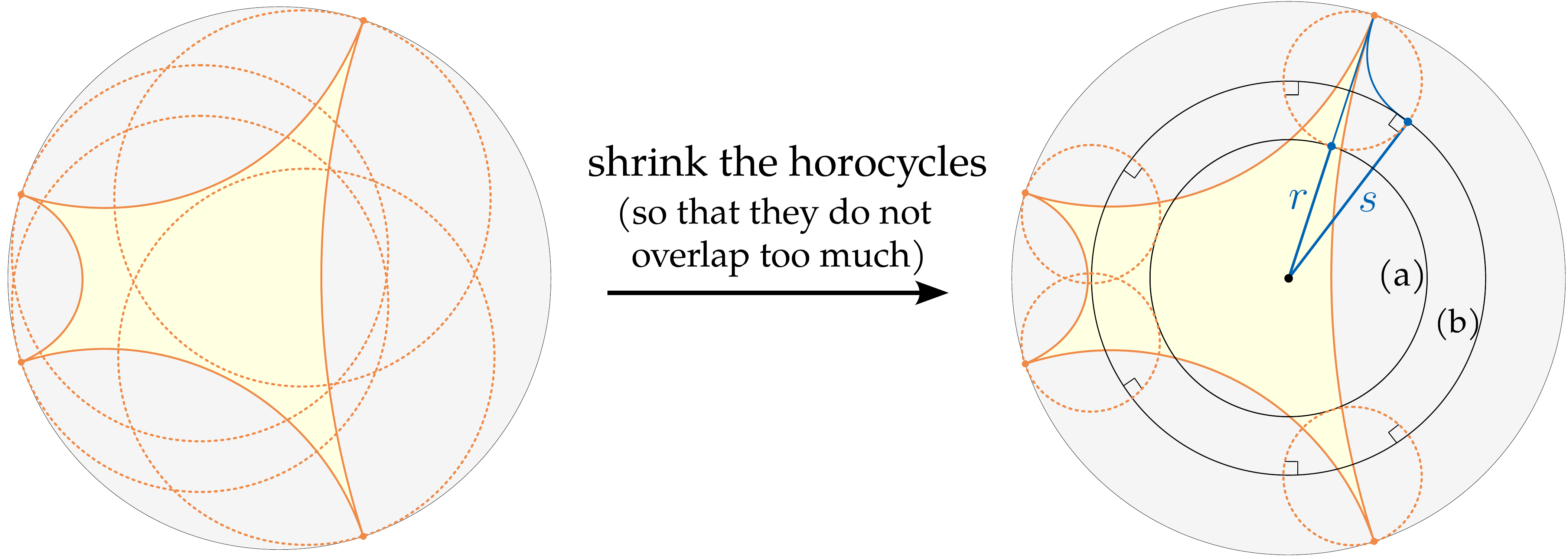}
	\caption{Conditions \ref{item_d1} and \ref{item_d2}.}
	\label{figure_cyclicdecorated}
\end{figure}
See Figure \ref{figure_cyclicdecorated}.
Conditions \ref{item_d1} and \ref{item_d2} are independent of the distance by which we shrink the horocycles, in the sense that if the desired circle exists for some shrinking distance $t$, then it also exists for any shrinking distance larger than $t$. Indeed, \ref{item_d1} can clearly be reformulated as:
\begin{enumerate}[label=(a')]
\item\label{item_d12} There exists a point in $\H^2$ whose signed distance\footnote{The signed distance from a point $p\in\H^2$ to a horocycle $H\subset\H^2$ is defined as the usual distance if $p$ lies outside of the horodisk enveloped by $H$, and otherwise is $-1$ times the usual distance.} to each of the horocycles is the same.
\end{enumerate}	 
Meanwhile, \ref{item_d2} is also equivalent to \ref{item_d12} by virtue of the fact that the lengths $r$ and $s$ in Figure \ref{figure_cyclicdecorated} are related by $r=\log \cosh s$.

\vspace{0.2cm}

Given adjacent cyclic decorated ideal polygons $\Delta'_1$ and $\Delta'_2$ which are \emph{compatible} in the sense that their horocycles coincide at both common vertices, the \emph{Delaunay condition} for $\Delta'_1$ and $\Delta'_2$ refers to either of the following equivalent conditions:
\begin{enumerate}[label=(\Alph*)]
    \item\label{item_dDelaunay1} After shrinking the horocycles of both $\Delta'_1$ and $\Delta'_2$ by the same distance if necessary, the circle tangent to the horocycles of $\Delta_1'$ as in Condition \ref{item_d1} is disjoint or tangent to each horocycle of $\Delta_2'$.
	\item\label{item_dDelaunay2} After shrinking the horocycles of both $\Delta'_1$ and $\Delta'_2$ by the same distance if necessary, the circle orthognoal to the horocycles of $\Delta_1'$ as in Condition \ref{item_d2} does not intersect any horocycle of $\Delta_2'$ at an angle larger than $\frac{\pi}{2}$.
	\item\label{item_dDelaunay3} For $i=1,2$, consider the planar polygon $\Delta''_i$ in $F$ whose vertices correspond to the horocycles of $\Delta_i'$ (cf.\ Condition \ref{item_d3}). Then $\Delta''_1$ and $\Delta''_2$ form a convex configuration (in the sense that  $\Delta''_1\cup\Delta''_2\subset\R^{2,1}$ is the graph of a convex function $x_3=f(x_1,x_2)$).
\end{enumerate}
\begin{figure}[h]
	\centering
	\includegraphics[width=16.4cm]{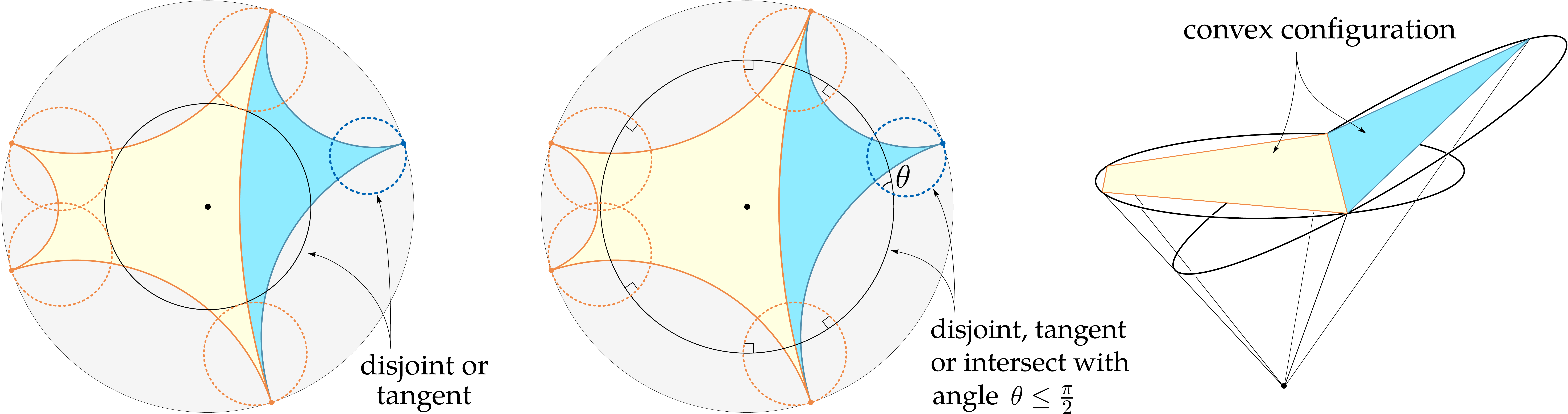}
	\caption{Conditions \ref{item_dDelaunay1}, \ref{item_dDelaunay2} and \ref{item_dDelaunay3}.}
	\label{figure_decoratedDelauny}
\end{figure}
See Figure \ref{figure_decoratedDelauny}. Again, Conditions \ref{item_dDelaunay1} and \ref{item_dDelaunay2} are independent of the shrinking distance and can be reformulated as
\begin{enumerate}[label=(A')]
	\item\label{item_dDelaunay12} Let $o_1'\in\H^2$ and $r_1'\in\R$ be the point and the signed distance in Condition \ref{item_d12} for $\Delta_1'$. Then the signed distance from $o_1'$ to any horocycle of $\Delta_2'$ is greater than or equal to $r_1'$.
\end{enumerate}	
We will take \ref{item_dDelaunay12} as our working definition and prove its equivalence with \ref{item_dDelaunay3} in Lemma \ref{lemma_FDelaunay}.

The \emph{strict} Delaunay condition for $\Delta'_1$ and $\Delta'_2$ are defined by the strict versions of the above conditions, whose formulations are obvious.
Using these conditions, we define Delaunay triangulation and decomposition for a decorated hyperbolic metric $(\ve{h},c)\in\dT_{g,n}$ in the same way as for cone-surfaces (Definition \ref{def_Delaunay}). Then a close analogue of Theorem \ref{thm_Delaunay} holds, namely every $(\ve{h},c)$ admits a unique Delaunay decomposition whose subdivisions are exactly the Delaunay triangulations.

\begin{remark}
Similarly as in the setting of hyperbolic cone-surfaces (see Remark \ref{remark_loose}), we could have loosened the definition of cyclic decorated ideal polygons, and hence the Delaunay condition, by allowing curves of constant geodesic curvature in place of circles in Conditions \ref{item_d1} and \ref{item_d2}, or equivalently, by removing the requirement ``$A\cap\pa F$ is an ellipse'' in \ref{item_d3}. This would make every decorated ideal triangle cyclic. Again, it does not change anything because any face of a Delaunay decomposition or triangulation under this looser definition is actually cyclic in our sense.
\end{remark}

Letting $(\ve{h},c)$ vary, we get the \emph{Penner decomposition} of the decorated Teichm\"uller space:
$$
\dT_{g,n}=\bigcup_{\T}\dT^\T_{g,n}~,
$$
where $\T$ runs through all homotopy classes of topological triangulations of $S_{g,n}$ with vertices at the marked points, and the \emph{Penner cell} $\dT^\T_{g,n}$ is the set of all $(\ve{h},c)\in\dT_{g,n}$ which admits a Delaunay triangulation of topological type $\T$. It is shown by Penner \cite{Penner} that the $\dT^\T_{g,n}$'s are the top dimensional closed cells of a cell decomposition of $\dT_{g,n}$, and by Akiyoshi \cite{Akiyoshi} that every fiber of the projection $\Pi:\dT_{g,n}\to\TT_{g,n}$ only meets finitely many such cells.

Finally, given a topological triangulation $\T$ with edges denoted by $\Edge_\T$, we recall that the \emph{Penner coordinates} $(\ell'_e)_{e\in\Edge_\T}:\dT_{g,n}\overset\sim\to\R^{\Edge_\T}$ are defined as follows. For any $(\ve{h},c)\in\dT_{g,n}$, we use the metric $\ve{h}$ to develop each edge $e\in\Edge_\T$ as a geodesic $l$ in $\H^2$ and develop the decorating horocycles $c$ at the two ends of $e$ as horocycle in $\H^2$ centered at the two ends of $l$. Then $\ell'_e(\ve{h},c)\in \R$ is defined to be the signed distance between the last two horocycles\footnote{About the notation: throughout the paper, we let $\ell$ represent the length of a line segment (usually a side of a polygon) in $\H^2$, $\E^2$ or $\SS^2$, which is  always positive; whereas the notation $\ell'$ (not derivative) is for the signed distance between two horocycles, which is positive/zero/negative when the horocyles are disjoint/tangent/intersecting (see e.g.\ \cite[Fig.\ 2]{Springborn}).}. A fundamental property is that two points of $\dT_{g,n}$ with Penner coordinates $(\ell'_e)$ and $(\wh{\ell}'_e)$ are in the same fiber of the projection $\Pi$ if and only if there exists $(u_i)\in\R^n$ such that 
$\hat{\ell}'_e=\ell'_e+u_i+u_j$ for all $i,j\in\{1,\cdots,n\}$ and all edge $e\in\Edge_\T$ joining the $i$th and $j$th marked points (see \cite[Prop.\ 3.2]{Springborn}).

\begin{proof}[Proof of Theorem \ref{thm_BPS}]
\textbf{Construction of $\wt{\Phi}_K$.} 
Given a convex cyclic polygon $\Delta$ in $\H^2$, $\E^2$ or $\SS^2$, the face $\Delta'$ of the polyhedron $P(\Delta)\subset\H^3$, developed as an ideal polygon in $\H^2$, has a natural decoration defined as follows: Consider the horosphere in $\H^3$ centered at each vertex of $\Delta'$ and furthermore
\begin{itemize}
	\item tangent to the horosphere or plane $L$ which contains $\Delta$ in the case of $\E^2$ or $\H^2$,
	\item passing through the point $o$ whose unit tangent sphere contains $\Delta$ in the case of $\SS^2$.
\end{itemize}
The intersection of this horosphere with the plane $L'\subset\H^3$ containing $\Delta'$ is a horocycle in $L'\cong\H^2$. These horocycles form a decoration of $\Delta'$, which we call the \emph{canonical} decoration. See Figure \ref{figure_canonicaldecoration}. 
\begin{figure}[h]
	\centering
	\includegraphics[width=11cm]{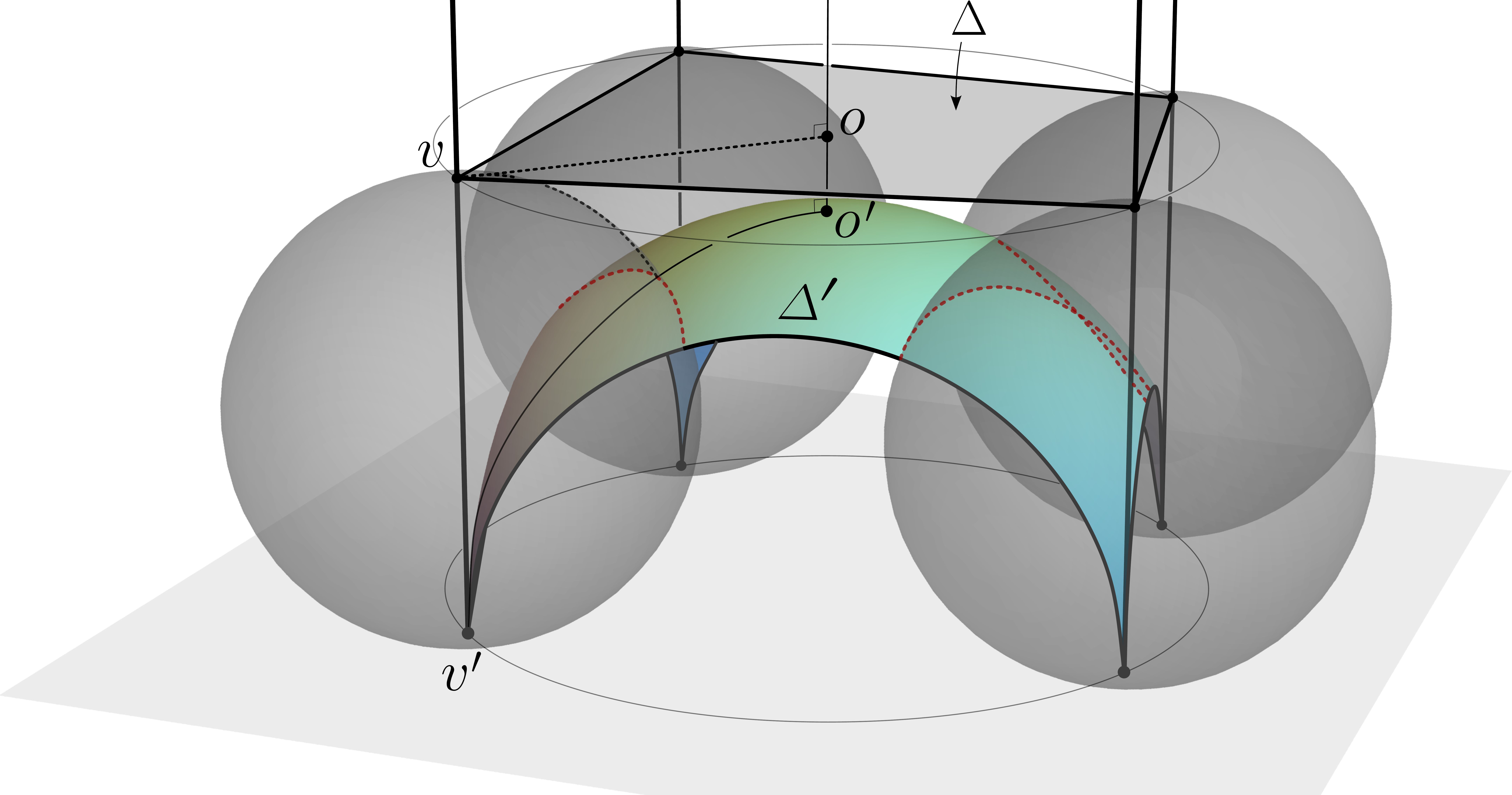}
	\caption{Canonical decoration of $\Delta'$ (the red horocycles) for a cyclic quadrilateral $\Delta\subset\E^2$.}
	\label{figure_canonicaldecoration}
\end{figure}
This gives $\Delta'$ the intrinsic structure of a cyclic decorated ideal polygon, where the cyclicity can be seen as follows: $L'$ has a distinguished point $o'$, namely the point closest to the above $L$ or $o$. Using trigonometric identities with horocycles (see \cite[Appendix A]{Guo-Luo_rigidity_II}; for example,  in the case of Figure \ref{figure_canonicaldecoration}, we apply the last identity in loc.\ cit.\  p.1306 with $\theta_3=\frac{\pi}{2}$, $\ell_3=0$ to the triangle $o'v'\infty$ together with the two horocycle arcs indicated in the figure), we see that the signed distance from $o'$ to the horosphere at every vertex of $\Delta'$ has the same value $r'$, which is related to the circumradius $r$ of $\Delta$ by
$$
e^{r'}=
\begin{cases}
\sin r &\text{if $\Delta\subset\SS^2$},\\
2r &\text{if $\Delta\subset\E^2$},\\
\sinh r &\text{if $\Delta\subset\H^2$}.
\end{cases}
$$
Thus, $\Delta'$ satisfies Condition \ref{item_d12}, as required.

Given adjacent convex cyclic polygons $\Delta_1,\Delta_2$, if we develop the faces $\Delta'_1,\Delta'_2$ of $P(\Delta_1)$, $P(\Delta_2)$ as adjacent ideal polygon in $\H^2$, then their canonical decorations are  clearly compatible. Therefore, for any $\ve{d}$ in $\M^{-1}_{g,n}$, $\M^0_{g,n}$ or $\M^{1,*}_{g,n}$, viewing $P(\ve{d})$ as the union $\bigcup_{\Delta\in\Face_\T}P(\Delta)$ where $\T$ is the Delaunay decomposition, we conclude that 
the canonical decorations on the $\Delta'$'s fit together to form a decoration of the punctured hyperbolic surface $\pa P(\ve{d})=\bigcup_{\Delta\in\Face_\T}\Delta'$. Pulling this decoration back to $S_{g,n}^\circ$ via the diffeomorphism $(S^\circ_{g,n},\ve{d})\to\pa P(\ve{d})$, we obtain a decoration $c$ for $\ve{h}=\Phi_K(\ve{d})\in\TT_{g,n}$. This defines the natural lift $\wt{\Phi}_K$ of $\Phi_K$.

\vspace{0.2cm}

\textbf{Injectivity.} We shall show that $\wt{\Phi}_K$ is injective, namely that $\ve{d}$ is completely determined by $\wt{\Phi}_K(\ve{d})=(\ve{h},c)$.
By Lemma \ref{lemma_Delaunay} below, the Delaunay decomposition of $(\ve{h},c)$ is also Delaunay for $\ve{d}$, so $(\ve{h},c)$ determines the combinatorics of the Delaunay decomposition of $\ve{d}$. Meanwhile, by Lemma \ref{lemma_trigonometry} (applied to the side faces of the polyhedron $P(\Delta)$ for each face $\Delta$ in the Delaunay decomposition of $(\ve{h},c)$), the Penner coordinates of $(\ve{h},c)$ also determine the edge lengths of the Delaunay decomposition of $\ve{d}$. It is a basic fact that any convex cyclic polygon in $\H^2$, $\E^2$ or $\SS^2$ is determined by its edge lengths (see e.g.\ \cite{Guo-Sonmez}). Therefore, we conclude that $(\ve{h},c)$ completely determines $\ve{d}$, as required.

\vspace{0.2cm}

\textbf{Characterization of discrete conformality.} 
The definition of discrete conformality can be reformulated as follows: two metrics in $\M^{-1}_{g,n}$, $\M^0_{g,n}$ or $\M^{1,*}_{g,n}$ are discretely conformal if and only if they are respectively the first and last members of a finite sequence in the same space, such that any adjacent members $\ve{d},\wh{\ve{d}}$ of the sequence
satisfy the following condition \autoref{eqn_proofBPS} for some topological triangulation $\T$ of $S_{g,n}$: 
\begin{equation}\label{eqn_proofBPS}\tag{$\star_{\scalebox{0.5}{$\T$}}$}
	\parbox{13.2cm}{
		$\T$ can be realized as Delaunay triangulations for both $\ve{d}$ and $\wh{\ve{d}}$. Moreover, there exists $(u_i)\in\R^n$ such that 
		for all $i,j\in\{1,\cdots,n\}$ and all $e\in\Edge_\T$ joining the $i$th and $j$th marked points, the geodesic lengths $\ell_e$ and $\hat{\ell}_e$ of $e$ under $\ve{d}$ and $\wh{\ve{d}}$ are related by  
		$$
		\hspace{-1cm}
		\begin{cases}
			\sinh\tfrac{\hat\ell_e}{2}=e^{\frac{u_i+u_j}{2}}\sinh\tfrac{\ell_e}{2}&\text{ for $\M^{-1}_{g,n}$}\\[5pt]    
			\hat\ell_e=e^{\frac{u_i+u_j}{2}}\ell_e&\text{ for $\M^0_{g,n}$,}\\[5pt]
			\sin\tfrac{\hat\ell_e}{2}=e^{\frac{u_i+u_j}{2}}\sin\tfrac{\ell_e}{2}&\text{ for $\M^{1,*}_{g,n}$.}
		\end{cases}
		$$
	}
\end{equation}

By the above mentioned coincidence between Delaunay decompositions of $(\ve{h},c)$ and $\ve{d}$ which follows from Lemma \ref{lemma_Delaunay}, as well as the relation between Penner coordinates of $(\ve{h},c)$ and the edge lengths of $\ve{d}$ which follows from Lemma \ref{lemma_trigonometry}, we infer that $\ve{d}$ and $\wh{\ve{d}}$ satisfy condition \autoref{eqn_proofBPS} if and only if the decorated hyperbolic metrics $(\ve{h},c):=\wt{\Phi}_K(\ve{d})$ and $(\wh{\ve{h}},\wh{c}):=\wt{\Phi}_K(\wh{\ve{d}})$ are both in the Penner cell $\dT_{g,n}^\T$ and their Penner coordinates $(\ell'_e)_{e\in\Edge_\T}$ and $(\wh{\ell}'_e)_{e\in\Edge_\T}$ are related by $\hat{\ell}'_e=\ell'_e+u_i+u_j$ for some $(u_i)\in\R^n$.
But as explained in the paragraph preceding this proof, this means exactly that $(\ve{h},c)$ and $(\wh{\ve{h}},\wh{c})$ are in the same fiber of $\Pi$ (namely, $\ve{h}=\wh{\ve{h}}$). Therefore, condition \autoref{eqn_proofBPS} on $\ve{d}$ and $\wh{\ve{d}}$ is equivalent to the following condition on $(\ve{h},c)$ and $(\wh{\ve{h}},\wh{c})$:
\begin{equation}\label{eqn_proofBPS3}\tag{$\star'_{\scalebox{0.5}{$\T$}}$}
	\parbox{13.5cm}{
		$(\ve{h},c)$ and $(\wh{\ve{h}},\wh{c})$ are both in the intersection of the Penner cell $\dT^\T_{g,n}$ with some fiber of $\Pi$.
	}
\end{equation}
Meanwhile, by the aforementioned finiteness result of Akiyoshi, two points of $\dT_{g,n}$ are in the same fiber of $\Pi$ if and only if they are respectively the first and last members of a finite sequence in $\dT_{g,n}$ such that any adjacent members $(\ve{h},c),(\wh{\ve{h}},\wh{c})$ of the sequence
satisfy  \autoref{eqn_proofBPS3} for some $\T$. Therefore, we conclude that any $\ve{d}_1,\ve{d}_2\in\M^{-1}_{g,n}$, $\M^0_{g,n}$ or $\M^{1,*}_{g,n}$ are discretely conformal if and only if $\wt{\Phi}_K(\ve{d}_1),\wt{\Phi}_K(\ve{d}_1)\in\dT_{g,n}$ are in the same fiber of $\Pi$, as required.
\end{proof}
The next two lemmas played a key role in the proof.
\begin{lemma}\label{lemma_Delaunay}
	Given adjacent convex cyclic polygons $\Delta_1$ and $\Delta_2$ in $\H^2$, $\E^2$ or $\SS^2$, consider the ideal polygonal faces $\Delta'_1$ and $\Delta'_2$ of the polyhedra $P(\Delta_1),P(\Delta_2)\subset\H^3$. Then the following conditions are equivalent:
	\begin{enumerate}[label=(\roman*)]
		\item\label{item_Dlemma1}
		$\Delta_1$ and $\Delta_2$ satisfy the Delaunay condition;	
		\item\label{item_Dlemma2} 
		the dihedral angle of the polyhedron $P(\Delta_1)\cup P(\Delta_2)$ at the common edge of $\Delta'_1$ and $\Delta'_2$ is at most $\pi$;
		\item\label{item_Dlemma3}
		$\Delta_1'$ and $\Delta_2'$, after being developed as adjacent ideal triangles in $\H^2$ and endowed with their canonical decorations, satisfy the Delaunay condition \ref{item_dDelaunay12}.
	\end{enumerate}
	Moreover, the strict versions of these conditions (i.e.\ replacing ``Delaunay'' by ``strict Delaunay'' in \ref{item_Dlemma1} \ref{item_Dlemma3}, and ``at most'' by ``less than'' in \ref{item_Dlemma2}) are also equivalent.
\end{lemma}
\begin{proof}
\textbf{``\ref{item_Dlemma1} $\Leftrightarrow$ \ref{item_Dlemma2}''.} The circumdisk $D_i$ of $\Delta_i$   ($i=1,2$) is cut by the edge $I:=\Delta_1\cap\Delta_2$ into two pieces, and we let $A_i$ be the piece containing $\Delta_i$. The Delaunay condition \ref{item_Dlemma1} is equivalent to the condition that the interior angle $\alpha$ of the bigon $A_1\cup A_2$ (which is bounded by two circle arcs) at either endpoint of $I$ is greater than or equal to $\pi$. See Figure \ref{figure_Delaunay}. 

On the other hand, identifying $\H^2$, $\E^2$ or $\SS^2$ with a plane/horosphere $L$ or the unit tangent sphere of a point $o$ in $\H^3$ as in \S \ref{subsec_relationship polyhedral}, we consider its map to the conformal boundary $\pa_\infty\H^3\cong\CP^1$ defined by orthogonal rays issuing from each point of it.
	This map is conformal and sends the circle $\pa D_i$ to the boundary at infinity of the plane $L_i'\subset\H^3$ containing $\Delta'_i$. By the conformality, the dihedral angle in  \ref{item_Dlemma2} is equal to $2\pi-\alpha$. The equivalence ``\ref{item_Dlemma1} $\Leftrightarrow$ \ref{item_Dlemma2}'' follows.
	 
\textbf{``\ref{item_Dlemma2} $\Leftrightarrow$ \ref{item_Dlemma3}''.} The plane $L_i'\subset\H^3$ containing $\Delta'_i$ is cut by the geodesic $I':=L_1'\cap L_2'$ into two half-planes, and we let $A_i'$ (resp.\ $B_i'$)  be the half-plane containing (resp.\ not containing) $\Delta'_i$. Let $\sigma\in\Isom^+(\H^3)$ be the hyperbolic rotation about $I'$ sending $A_2'$ to $B_1'$. See Figure \ref{figure_dihedral}.
\begin{figure}[h]
	\centering
	\includegraphics[width=17cm]{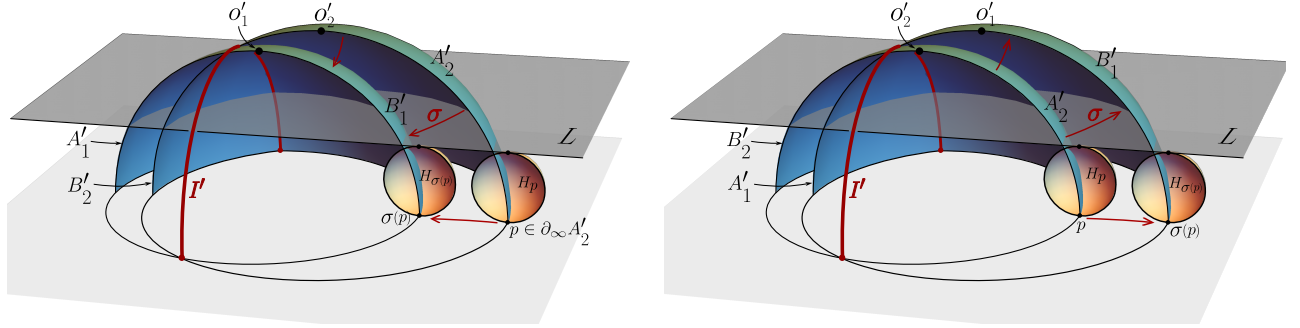}
	\caption{To the proof of ``\ref{item_Dlemma2} $\Leftrightarrow$ \ref{item_Dlemma3}'' when $\Delta_1,\Delta_2\subset\E^2$. On the left and right are the cases with dihedral angle greater and less than $\pi$, respectively.}
	\label{figure_dihedral}
\end{figure}

For each point $p$ on the boundary $\pa_\infty L_i'\subset\pa_\infty\H^3$ of $L_i'$ at infinity, consider the horosphere $H_p$ centered at $p$ and tangent to $L$ or passing through $o$. This defines a family of horospheres with the same signed distance $r_i'$ to the distinguished point $o_i'$ of $L_i'$, and the canonical decoration of $\Delta_i'$ is given by these horospheres (see the proof of Theorem \ref{thm_BPS}). 
The Delaunay condition \ref{item_dDelaunay12} is equivalent to the condition that $r_1'$ is less than or equal to the signed distance from $o_1'$ to any $\sigma(H_p)$ with $p\in \pa_\infty A_2'$. On the other hand, by looking at the plane orthogonal to $I'$ and containing $p$ at infinity as in Figure \ref{figure_dihedral}, we see that if the dihedral angle in \ref{item_Dlemma2} is $\leq\pi$ (resp.\ $\geq\pi$), then 
the horosphere  $\sigma(H_p)$ is contained in $H_{\sigma(p)}$ (resp.\ contains  $H_{\sigma(p)}$). The required equivalence follows.

Finally, the ``Moreover'' statement is proved by following the same arguments. We omit the details.
\end{proof}

\begin{figure}[h]
	\centering
	\includegraphics[width=13.3cm]{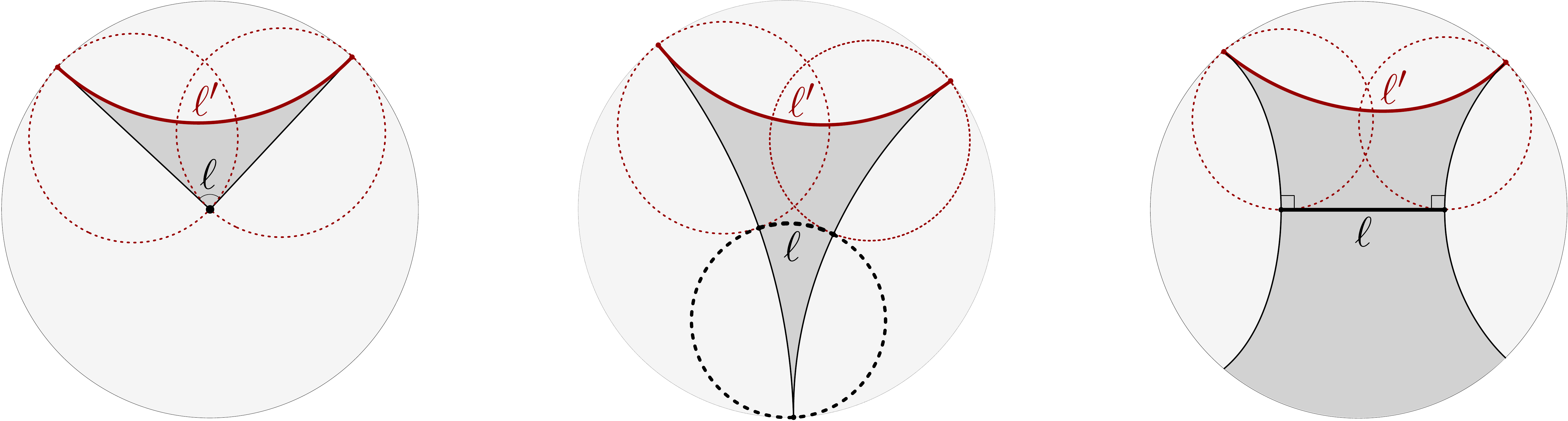}
	\caption{Configurations of geodesics (solid) and horocycles (dashed) in $\H^2$ exhibited by each side face of $P(\Delta)$, when $\Delta$ is in $\SS^2$, $\E^2$ and $\H^2$, respectively. 
		The signed distance $\ell'$ between the red horocycles is negative in these examples, but can also be positive in the latter two cases (if the boldfaced horocycle/geodesic is close enough to the red geodesic).
	}
	\label{figure_trigonometry}
\end{figure}

\begin{lemma}\label{lemma_trigonometry}
	For the three configurations in Figure \ref{figure_trigonometry} respectively, we have
	$$
	e^{\frac{\ell'}{2}}=
	\begin{cases}
		\sin\frac{\ell}{2}&\text{ (left),}\\ 
		\ell&\text{ (middle),}\\ 
		\sinh\frac{\ell}{2}&\text{ (right).}
	\end{cases}
	$$
\end{lemma}
\begin{proof}
	These are special cases of more general hyperbolic trigonometric identities with horocycles given e.g.\ in  \cite[Appendix A]{Guo-Luo_rigidity_II}: the three figures here are exactly the last three therein with $l_1=l_2=0$.
\end{proof}

\section{Epstein-Penner convex hulls in $\R^{2,1}$, $\SS^{2,1}$ and $\H^{2,1}$}\label{sec_3}
\subsection{Preliminaries on $\SS^{2,1}$ and $\H^{2,1}$}\label{subsec_preliminary}
In the projective space $\RP^3$, there are two projective equivalence classes of nondegenerate quadric surfaces,
namely  \emph{ellipsoids} and  \emph{doubly ruled quadrics}, which can be identified respectively with the sphere $\{x_1^2+x_2^2+x_3^2=1\}$ and the one-sheeted hyperboloid $\{x_1^2+x_2^2-x_3^2=1\}$ in a suitable affine chart $\R^3\subset\RP^3$.
By definition, the interior and exterior domains of an ellipsoid are the projective models of the hyperbolic space $\H^3$ and the \emph{de Sitter space} $\SS^{2,1}$, respectively, whereas the two connected components of $\RP^3$ minus a doubly ruled quadric are projectively equivalent to each other, and either of them is a projective model for the \emph{anti-de Sitter space} $\H^{2,1}$.

$\SS^{2,1}$ (resp.\ $\H^{2,1}$) carries a canonical Lorentzian metric $\ve{g}_{1}$ (resp.\ $\ve{g}_{-1}$) of constant sectional curvature $1$ (resp.\ $-1$), with the properties that the isometries are exactly the projective transformations of $\RP^3$ preserving $\SS^{2,1}$ (resp.\ $\H^{2,1}$) as a subset, and that the geodesic lines/planes are exactly the projective lines/planes.
Both metrics are constructed similarly as the hyperboloid construction for the Riemannian metric on $\H^3$:
Let $\R^{p,q}$ denote the vector space $\R^{p+q}$ equipped with the scalar product
$$
\langle \xi,\eta\rangle:=\xi_1\eta_1+\cdots+\xi_p\eta_p-\xi_{p+1}\eta_{p+1}-\cdots-\xi_{p+q}\eta_{p+q}~.
$$
Then the sets of null vectors in $\R^{3,1}$ and $\R^{2,2}$ projectivize to an ellipsoid and a  doubly ruled quadric, respectively, in $\RP^3=\P(\R^{3,1})=\P(\R^{2,2})$, whereas $\SS^{2,1}$ (resp.\ $\H^{2,1}$) is projectivized from the set of spacelike (resp.\ timelike) vectors in $\R^{3,1}$ (resp.\ $\R^{2,2}$), and hence identifies with the antipodal quotient of the hypersurface $\{\langle\xi,\xi\rangle=1\}$ (resp.\ $\{\langle\xi,\xi\rangle=-1\}$). Namely, we may write
$$
\SS^{2,1}:=\big\{[\xi]\in\P(\R^{3,1})\,\big|\,\langle\xi,\xi\rangle>0\big\}\cong\big\{\xi\in\R^{3,1}\,\big|\, \langle \xi,\xi\rangle=1\big\}/\pm\!1,
$$
$$
\H^{2,1}:=\big\{[\xi]\in\P(\R^{2,2})\,\big|\,\langle\xi,\xi\rangle<0\big\}\cong\big\{\xi\in\R^{2,2}\,\big|\, \langle \xi,\xi\rangle=-1\big\}/\pm\!1.
$$
The induced metrics on these hypersurfaces then define the metrics $\ve{g}_{\pm1}$ of $\SS^{2,1}$ and $\H^{2,1}$.

The Euclidean counterpart of the Lorentzian manifolds $\SS^{2,1}$ and $\H^{2,1}$ is just the Minkowski space $\R^{2,1}$. We will denote points in $\R^{2,1}$ by $x$, $y$ etc., the scalar product by $x\pro y:=x_1y_1+x_2y_2-x_3y_3$, and the Lorentzian metric by $\ve{g}_0:=\dx_1^2+\dx_2^2-\dx_3^2$.
Meanwhile, we also use the inclusion
$$
\R^{2,1}\hookrightarrow \RP^3=\P(\R^{3,1})=\P(\R^{2,2}),\ \ x\mapsto [\xi]=[x_1:x_2:1:x_3]
$$
to view $\R^{2,1}$ as an affine chart in $\RP^3$  (this choice of chart is nonstandard for $\SS^{2,1}$, cf.\ Remark \ref{remark_chart1}). The parts of $\SS^{2,1}$ and $\H^{2,1}$ in this chart are respectively the regions
$$
\Omega_1:=\SS^{2,1}\cap \R^{2,1}=\big\{x\in\R^{2,1}\,\big|\,x\pro x+1>0\big\},
$$
$$ 
 \Omega_{-1}:=\H^{2,1}\cap \R^{2,1}=\big\{x\in\R^{2,1}\,\big|\,x\pro x-1<0\big\}.
$$
See Figure \ref{figure_Sitter}. In other words, we have identified $\SS^{2,1}$ and $\H^{2,1}$, with a plane removed, as the regions $\Omega_1$ and $\Omega_{-1}$ in $\R^{2,1}$, respectively. By straightforward calculation, we get the following expressions for the Lorentizan metrics $\ve{g}_{\pm1}$ on these regions:
$$
\ve{g}_1=\frac{1}{(1+x\pro x)^2}
\scalebox{0.9}{$
\begin{pmatrix}
1+x_2^2-x_3^2&-x_1x_2&x_1x_3\\[0.1cm]
-x_1x_2&1+x_1^2-x_3^2&x_2x_3\\[0.1cm]
x_1x_3&x_2x_3&-1-x_1^2-x_2^2
\end{pmatrix}
$}\quad\text{for }x\in\Omega_1,
$$
\hspace{0.1cm}
$$
\ve{g}_{-1}=\frac{1}{(1-x\pro x)^2}
\scalebox{0.9}{$
\begin{pmatrix}
	1-x_2^2+x_3^2&x_1x_2&-x_1x_3\\[0.1cm]
	x_1x_2&1-x_1^2+x_3^2&-x_2x_3\\[0.1cm]
	-x_1x_3&-x_2x_3&-1+x_1^2+x_2^2
\end{pmatrix}
$}\quad\text{for }x\in\Omega_{-1},
$$
where a symmetric matrix-valued function $A(x)$ is understood as the metric $\transp\dx A(x)\dx$.

\begin{figure}[h]
	\centering
	\includegraphics[width=5.6cm]{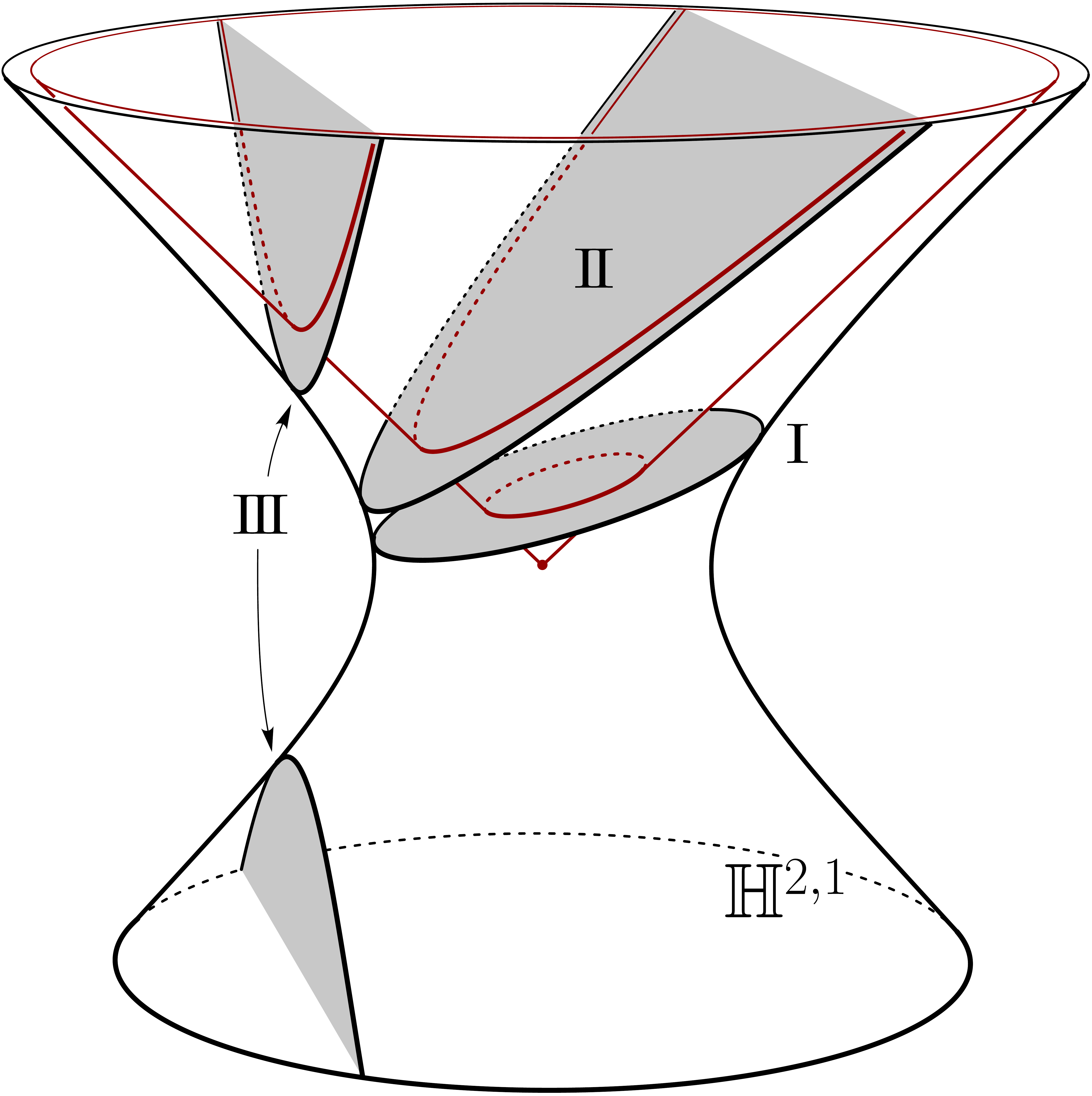}
	\caption{The three types of spacelike planes in $\H^{2,1}$ relative to the affine chart. A plane of type II (resp.\ III) meets the plane at infinity at a point (resp.\ along a line).
	}
	\label{figure_adsplanes}
\end{figure}

In particular, we view the origin $0\in\R^{2,1}$ also as a point of both $\SS^{2,1}$ and $\H^{2,1}$, at which the metrics $\ve{g}_{\pm1}$ both coincide with the Minkowski metric $\ve{g}_0$ of $\R^{2,1}$. By the following lemma, the stabilizer of this point in the isometry group of $\SS^{2,1}$ or $\H^{2,1}$ is always the group $\SO(2,1)$ of linear automorphisms of $\R^{2,1}$, and the action of this stabilizer also coincides with the linear one:
\begin{lemma}\label{lemma_linear}
The linear action of $\SO(2,1)$ on $\R^{2,1}$ preserves the regions $\Omega_{\pm1}$ and the Lorentzian metrics $\ve{g}_{\pm1}$ on them.
\end{lemma}
\begin{proof}
The linear action preserves $\Omega_{\pm1}$ because these regions are defined by the Minkowski scalar product ``$\pro$'', which is preserved by the action. To see that $\ve{g}_{\pm1}$ are also preserved, note that the linear action of any $A\in\SO(2,1)$ can also be viewed as a projective transformation of $\RP^4$ preserving the affine chart $\R^{2,1}\subset\RP^4$ and the subsets $\SS^{2,1},\H^{2,1}\subset\RP^4$ at the same time, but such a projective transformation preserves $\ve{g}_{\pm1}$.
\end{proof}
Finally, the (geodesic) spacelike  lines and planes in $\SS^{2,1}$ and $\H^{2,1}$ can be characterized as follows:
\begin{itemize}
	\item The spacelike lines (resp.\ planes) in $\SS^{2,1}$ are exactly the lines (resp.\ planes) in $\RP^3$  entirely contained in $\SS^{2,1}$. Every such plane is isometric to the projective plane $\RP^2$ with the round metric (i.e.\ the antipodal quotient of $\mathbb{S}^2$).
	\item 
	The spacelike lines (resp.\ planes) in $\H^{2,1}$ are exactly the intersection parts of $\H^{2,1}$ with projective lines (resp.\ planes) in $\RP^3$ which are open segments (resp.\ ellipses). Every such plane is isometric to $\H^2$, but can be in three types of relative positions with respect to the affine chart, see Figure \ref{figure_adsplanes}.
\end{itemize}

\begin{remark}\label{remark_chart1}
A more commonly used affine chart of $\RP^3$ for $\SS^{2,1}$ is given by the inclusion
$$
\R^3\hookrightarrow\RP^3=\P(\R^{3,1}),\quad x\mapsto [\xi]=[x_1:x_2:x_3:1].
$$
The part of $\SS^{2,1}$ in this $\R^3$ is the ball-complement
$$
\wt{\Omega}_1=\big\{x\in\R^3\,\big|\, x_1^2+x_2^2+x_3^2>1\big\}.
$$
The projective involution $[\xi_1:\xi_2:\xi_3:\xi_4]\mapsto[\xi_1:\xi_2:\xi_4:\xi_3]$, which intertwines the two charts, is expressed in the coordinates of either chart as 
\begin{equation}\label{eqn_involution}
x\mapsto y=\Big(\frac{x_1}{x_3},\frac{x_2}{x_3},\frac{1}{x_3}\Big).
\end{equation}
This can be understood as a transformation both from the region $\Omega_1$ to $\wt{\Omega}_1$ and back. On the other hand, the part of $\H^{2,1}$ in this chart is still $\Omega_{-1}$, or in other words, the transformation \autoref{eqn_involution} preserves $\H^{2,1}$. See also Remark \ref{remark_chart2} below.
\end{remark}

\subsection{The cone $F$ and polygons in it}\label{subsec_F}
In this subsection, we collect some key facts about the cone
$$
F:=\big\{x\in\R^{2,1}\,\big|\,x\pro x\leq0,\ x_3\geq0\big\}
$$
and planar polygons in $F$ spanned by vertices on the boundary $\pa F=\{x\pro x=0,\,x_3\geq0\}$.

Since $F$ is contained in the region $\Omega_{-1}=\H^{2,1}\cap\R^{2,1}$ defined above, we view it as a cone not only in the Lorentizan manifold $\R^{2,1}$, but also in $\H^{2,1}$. Meanwhile, its subset
$$
F^*:=\big\{x\in\R^{2,1}\,\big|\,-1<x\pro x\leq0,\ x_3\geq0\big\}
$$
is contained in $\Omega_{1}$ and hence in $\SS^{2,1}$.
When viewed as a subset of $\R^{2,1}$, $\H^{2,1}$ or $\SS^{2,1}$, $F$ or $F^*$ is endowed with the Lorentizan metric $\ve{g}_0$, $\ve{g}_{-1}$ or $\ve{g}_1$, respectively.
In any case,  $F$ or $F^*$ is the \emph{future cone} of the origin $0$ in the respective Lorentzian manifold, namely the union of all \emph{future oriented causal paths} issuing from $0$ (see e.g.\ \cite{Bonsante-Seppi} for definition). 


It follows from Lemma \ref{lemma_linear} that the stabilizer of $F$ or $F^*$ in the isometry group of $\R^{2,1}$, $\H^{2,1}$ or $\SS^{2,1}$ is always the identity component $\SO_0(2,1)\cong\PSL(2,\R)$ of $\SO(2,1)$, and moreover its action is always the linear one in all three cases. We have the following basic fact about the action on the boundary $\pa F$:
\begin{lemma}\label{lemma_actionpair}
	For any $t>0$, the diagonal action of $\SO_0(2,1)$ on the subset $G_t:=\big\{(x,y)\,\big|\, x\pro y=-t\big\}$ of $\pa F\times\pa F$ is free and transitive.
\end{lemma}
Note that we have $x\pro y\leq 0$ for all $(x,y)\in \pa F\times \pa F$, with equality only when $x=0$ or $y=0$. Therefore, the sets $(G_t)_{t>0}$ form a foliation of $(\pa F\setminus\{0\})\times (\pa F\setminus\{0\})$.
\begin{proof}
Under the duality between points of $\pa F\setminus\{0\}$ and horocycles in $\H^2$ shown in Figure \ref{figure_duality}, it is a simple fact that the Minkowski scalar product of $x,y\in \pa F\setminus\{0\}$ is related to the signed distance $\ell'$ between the corresponding horocycles $x^*,y^*\subset\H^2$ by 
$x\pro y=-2e^{\ell'}$ (see \cite[Lemma 2.1]{Penner}). Therefore, in the space of ordered pairs of horocycles, $G_t$ corresponds to the subset formed by pairs with fixed signed distance $\log\frac{t}{2}$. It is an elementary fact in hyperbolic geometry that the action of $\Isom^+(\H^2)\cong\SO_0(2,1)$ on such a subset is free and transitive. Using the duality (which respects the $\SO_0(2,1)$-actions) to transfer back, we conclude that so is the action on $G_t$.
\end{proof}

\begin{figure}[h]
	\centering
	\includegraphics[width=15.5cm]{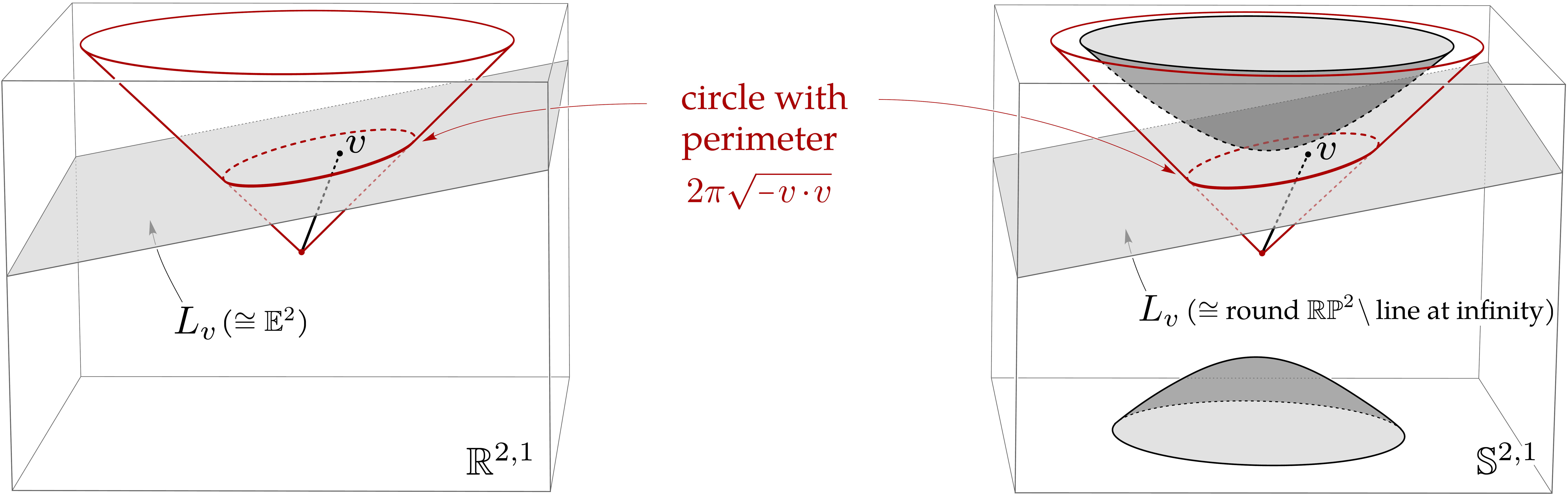}
	\caption{Lemma \ref{lemma_spacelike} in the cases of $\R^{2,1}$ and $\SS^{2,1}$ (see Figure \ref{figure_adsplanes} for $\H^{2,1}$).
	}
	\label{figure_planes}
\end{figure}

We often consider affine planes in $\R^{2,1}$ intersecting $\pa F$ along ellipses. They are exactly planes of the form
$$
L_v:=v+v^\perp=\big\{v+x\,\big|\,x\in\R^{2,1},\ x\pro v=0\big\}
$$
with $v$ in the interior of $F$. The induced metric on such a plane can be understood as follows:
\begin{lemma}\label{lemma_spacelike}
For any point $v$ in the interior of $F$, the affine plane $L_v$ is spacelike in $\R^{2,1}$, whereas $L_v\cap \Omega_{-1}$ is spacelike in $\H^{2,1}$. If $L_v$ is contained in $\Omega_1$, or equivalently if $v\pro v>-1$, then $L_v$ is spacelike in $\SS^{2,1}$ as well. Moreover, 
\begin{itemize}
\item
With the metric induced from the Lorentzian metric $\ve{g}_0$ of $\R^{2,1}$, $L_v$ is isometric to $\E^2$, and $L_v\cap F$ is the disk therein centered at $v$ with radius $r=\sqrt{-v\pro v}$.
\item
With the metric induced from the Lorentzian metric $\ve{g}_{-1}$ of  $\H^{2,1}$, $L_v\cap\Omega_{-1}$ is isometric to $\H^2$, and $L_v\cap F$ is the disk therein centered at $v$ with radius $r$ given by $\sinh  r=\sqrt{-v\pro v}$.
\item
Suppose $v\pro v>-1$. Then with the metric induced from the Lorentzian metric $\ve{g}_1$ of  $\SS^{2,1}$, $L_v$ is isometric to an open hemisphere in $\SS^2$ (or in other words, the round $\RP^2$ with a line removed), and $L_v\cap F$ is the disk therein centered at $v$ with radius $r$ given by $\sin  r=\sqrt{-v\pro v}$.
\end{itemize}
In any case, the perimeter of the circle $L_v\cap\pa F$ is $2\pi\sqrt{-v\pro v}$ (see Figure \ref{figure_planes}).
\end{lemma}
\begin{proof}
The spacelikeness and the description of $L_v$ and $L_v\cap\Omega_{-1}$ as $\E^2$, $\H^2$ and hemisphere follow from the discussions preceding Remark \ref{remark_chart1}. Since the perimeter of a circle with radius $r$ in $\E^2$, $\H^2$ and $\SS^2$ is equal to $2\pi r$, $2\pi \sinh r$ and $2\pi\sin r$, respectively, the desired formulas for radius and the one for perimeter are equivalent. Therefore, it only remains to show the formula for perimeter. To this end, we transform $v$ into the vertical vector $v_0=(0,0,\sqrt{-v\pro v})$ with an element of $\SO_0(2,1)$. By Lemma \ref{lemma_linear}, the perimeter of $L_v\cap\pa F$ is the same as that of the horizontal circle $L_{v_0}\cap\pa F$, while it is straightforward to calculate the latter and get the required value $2\pi\sqrt{-v\pro v}$ by using the expressions of $\ve{g}_{\pm1}$ in \S \ref{subsec_preliminary} and the rotational symmetry.
\end{proof}

Planar polygons in $F$ circumscribed by ellipses of the form $L_v\cap\pa F$ play a key role in the Epstein-Penner construction. 
We always denote such a polygon by $\Delta''$ (while the notation $\Delta$ is reserved for compact polygons in $\E^2$, $\H^2$ or $\SS^2$, and $\Delta'$ is for decorated ideal polygons in $\H^2$). Lemma \ref{lemma_spacelike} implies that any $\Delta''$, endowed with the metric induced from $\ve{g}_0$, $\ve{g}_{-1}$ or $\ve{g}_1$ (the last makes sense only when $\Delta''\subset F^*$), is intrinsically a convex cyclic polygon in $\E^2$, $\H^2$ or $\SS^2$ in the sense of \S \ref{subsec_21}. 
Conversely, the following lemma says that every convex cyclic polygon can be realized as some $\Delta''$, which is unique up to the $\SO_0(2,1)$-action\footnote{The uniqueness here is understood in the sense of \emph{marked} polygons. Namely, polygons should be considered as endowed with a labeling of vertices by $1,\cdots,n$ consecutively, while isomorphisms of polygons are required to respect the labels.}: 
\begin{lemma}\label{lemma_realizetriangle}
	Let $\Delta$ be a convex cyclic polygon in $\E^2$ or $\H^2$ (resp.\ $\SS^2$). Then there exists, up to the action of $\SO_0(2,1)$, a unique planar polygon $\Delta''$ in $F$ (resp.\ $F^*$) with vertices in $\pa F$, such that when endowed with the metric induced from $\ve{g}_0$ or $\ve{g}_{-1}$ (resp.\ $\ve{g}_1$), $\Delta''$ is isometric to $\Delta$. The affine plane containing $\Delta''$ intersects $\pa F$ along an ellipse.
\end{lemma}
\begin{proof}
\textbf{Existence.} Given $\Delta$, consider its circumdisk $D$ in $\E^2$, $\H^2$ or $\SS^2$. Using Lemma \ref{lemma_spacelike}, we find an affine plane $L_v\subset\R^{2,1}$ such that $L_v\cap F$ is isometric to $D$ (when endowed with the metric $\ve{g}_0$, $\ve{g}_{-1}$ or $\ve{g}_{1}$, respectively). We then obtain the required $\Delta''$ as the image of $\Delta$ under an isometry $D\overset\sim\to L_v\cap F$.

\textbf{The last statement.} We shall show that if a planar polygon $\Delta''\subset F$ has vertices in $\pa F\setminus\{0\}$ and is isometric to some convex cyclic polygon $\Delta\subset\E^2$, $\H^2$ or $\SS^2$ (when endowed with $\ve{g}_0$, $\ve{g}_{-1}$ or $\ve{g}_{1}$, respectively, assuming $\Delta''\subset F^*$ in the case of $\ve{g}_1$), then the plane $L$ containing $\Delta''$ must intersect $\pa F$ along an ellipse. 

Suppose otherwise. Then $L\cap\pa F$ is a parabola or hyperbola. In the Euclidean or spherical case, by the discussion preceding Remark \ref{remark_chart1}, $L$ is non-spacelike in the sense that any open subset (in particular the interior of $\Delta''$) contains lightlike and timelike line segments, a contradiction. In the hyperbolic case, $L\cap\Omega_{-1}$ is either non-spacelike in the same way, which leads to contradiction, or can still be spacelike. In the latter situation, if we add the points at infinity to $L\cap\Omega_{-1}$ and view it as a copy of $\H^2$ in $\H^{2,1}$, then the parabola or hyperbola $L\cap \pa F$ is a horocycle or an equidistance curve to a geodesic in this $\H^2$ and circumscribes $\Delta''$ (see Figure \ref{figure_adsplanes}), contradicting the assumption that $\Delta$ is cyclic.

\textbf{Uniqueness up to $\SO_0(2,1)$.} Suppose $\Delta''_1,\Delta''_2\subset F$ are both isometric to a given $\Delta$. By the last statement that we just proved, there are vectors $v_1,v_2$ in the interior of $F$ such that $\Delta_i''$ is circumscribed by the ellipse $L_{v_i}\cap\pa F$ ($i=1,2$). We have $v_1\pro v_1=v_2\pro v_2$ because by Lemma \ref{lemma_spacelike} this value only depends on the radius of the circumdisk of $\Delta$ . It follows that $v_1$ can be brought to $v_2$  by some $A\in\SO_0(2,1)$, or equivalently, $L_{v_1}$ can be brought to $L_{v_2}$ by $A$. Composing $A$ with a rotation of $L_{v_2}$ about $v_2$, we can bring $\Delta''_1$ to $\Delta''_2$, as required. 
\end{proof}

Besides the intrinsic structure of a convex cyclic polygon $\Delta$ in $\E^2$, $\H^2$ or $\SS^2$, the above $\Delta''$ also carries the data of a decorated ideal polygon $\Delta'$ in $\H^2$. In fact, the radial projection sends $\Delta''$ homeomorphically to an ideal polygon in the hyperboloid $\{x\pro x=-1\}\cong\H^2$, while the location of the vertices of $\Delta''$ defines a decoration via the duality in Figure \ref{figure_duality}. We will call this $\Delta'$ the decorated ideal polygon \emph{underlying} $\Delta''$. Note that $\Delta'$ is cyclic in the sense of \S \ref{subsec_relationship decorated}.

\begin{figure}[h]
	\centering
	\includegraphics[width=13.5cm]{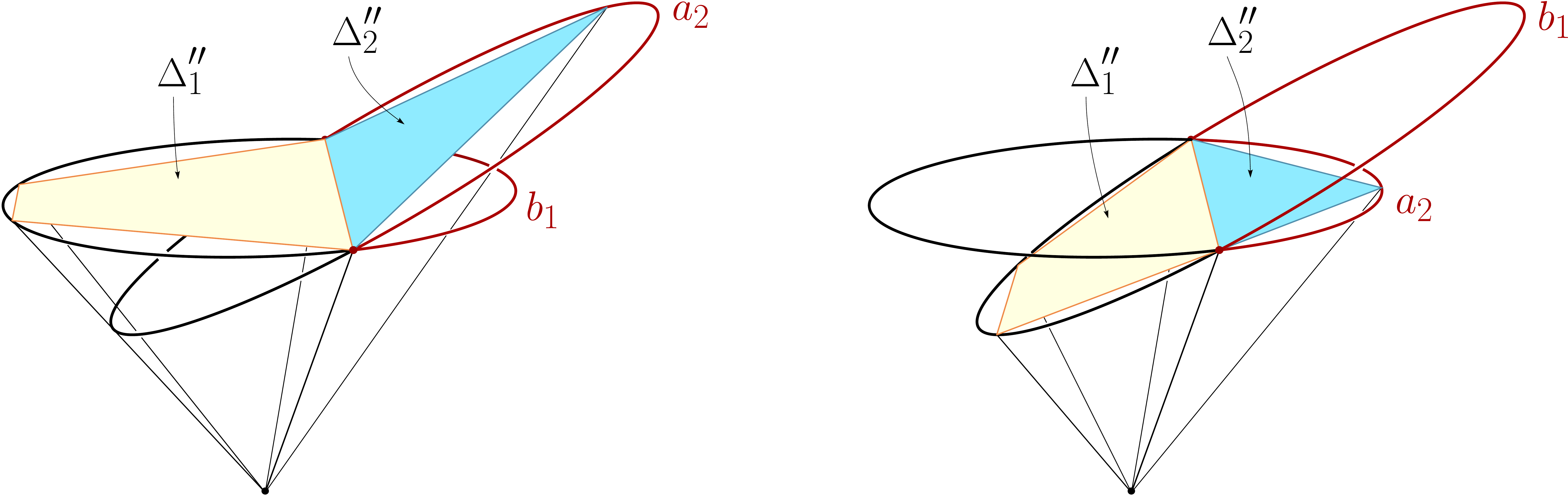}
	\caption{$\Delta_1'',\Delta_2''\subset F$ in convex configuration (left) and concave configuration (right).
	}
	\label{figure_convexconfiguration}
\end{figure}

In essence, this paper is based on the fact that the link between $\Delta$ and $\Delta'$ through the polygon $\Delta''\subset F$ coincides with the one in \S \ref{sec_2} through the polyhedron $P(\Delta)\subset\H^3$ (see also Remark \ref{remark_hilbert}), modulo a dilation in the Euclidean case. Corollary \ref{coro_maincoro} below gives a precise statement. In particular, the next two lemmas about the former link are the counterpart of Lemmas \ref{lemma_Delaunay} and \ref{lemma_trigonometry} about the latter.
\begin{lemma}\label{lemma_FDelaunay}
Let $\Delta_1$ and $\Delta_2$ be adjacent convex cyclic polygons in $\H^2$, $\E^2$ or $\SS^2$. Consider adjacent planar polygons $\Delta''_1,\Delta''_2\subset F$ isometric to them (given by Lemmas \ref{lemma_actionpair} and \ref{lemma_realizetriangle})
and the decorated ideal polygon $\Delta'_1$ and $\Delta'_2$ underlying $\Delta_1''$ and $\Delta_2''$.
	Then the following conditions are equivalence to each other:
	\begin{enumerate}[label=(\roman*)]
		\item\label{item_Flemma1}  
		$\Delta_1$ and $\Delta_2$ satisfy the Delaunay condition;	
		\item\label{item_Flemma2}  
		$\Delta_1''$ and $\Delta_2''$ form a convex configuration (see \S \ref{subsec_relationship decorated} and Figure \ref{figure_convexconfiguration});
		\item\label{item_Flemma3} 
		$\Delta_1'$ and $\Delta_2'$, after being developed in $\H^2$ as adjacent cyclic decorated ideal polygons, satisfy the Delaunay condition \ref{item_dDelaunay12} in \S \ref{subsec_relationship decorated}.
	\end{enumerate}
	Moreover, the strict versions of these conditions (i.e.\ replacing ``Delaunay'' by ``strict Delaunay'' in \ref{item_Flemma1} \ref{item_Flemma3}, and ``convex'' by ``strictly convex'' in \ref{item_Flemma2}) are also equivalent.
\end{lemma}
The equivalence ``\ref{item_Flemma2}$\Leftrightarrow$\ref{item_Flemma3}'' signifies that the formulation  \ref{item_dDelaunay3} of the Delaunay condition for ideal polygons is equivalent to the other formulations in \S \ref{subsec_relationship decorated}.
\begin{proof}
\textbf{``\ref{item_Flemma1}$\Leftrightarrow$\ref{item_Flemma2}''.} The circumcircle of $\Delta_i$ ($i=1,2$) is cut by the common vertices of $\Delta_1$ and $\Delta_2$ into two arcs, and we let $a_i$ (resp.\ $b_i$) be the length of the arc containing (resp.\ not containing) the other vertices of $\Delta_i$. Then the Delaunay condition \ref{item_Flemma1} is equivalent to $a_2\geq b_1$, see Figure \ref{figure_Delaunay}.

On the other hand, $a_2$ and $b_1$ are also the lengths of the ellipse arcs on $\pa F$ indicated in Figure \ref{figure_convexconfiguration} (under the metric $\ve{g}_{-1}$, $\ve{g}_0$ or $\ve{g}_1$ when $\Delta_1$ and $\Delta_2$ are in $\H^2$, $\E^2$ or $\SS^2$, respectively). Therefore, to establish the equivalence between \ref{item_Flemma1} and \ref{item_Flemma2}, it suffices to show that for 
any distinct $x,y\in\pa F\setminus\{0\}$, the one-parameter family of ellipse arcs passing through $x$ and $y$, as shown in Figure \ref{figure_arcs}, have strictly increasing lengths.
\begin{figure}[h]
	\centering
	\includegraphics[width=7.6cm]{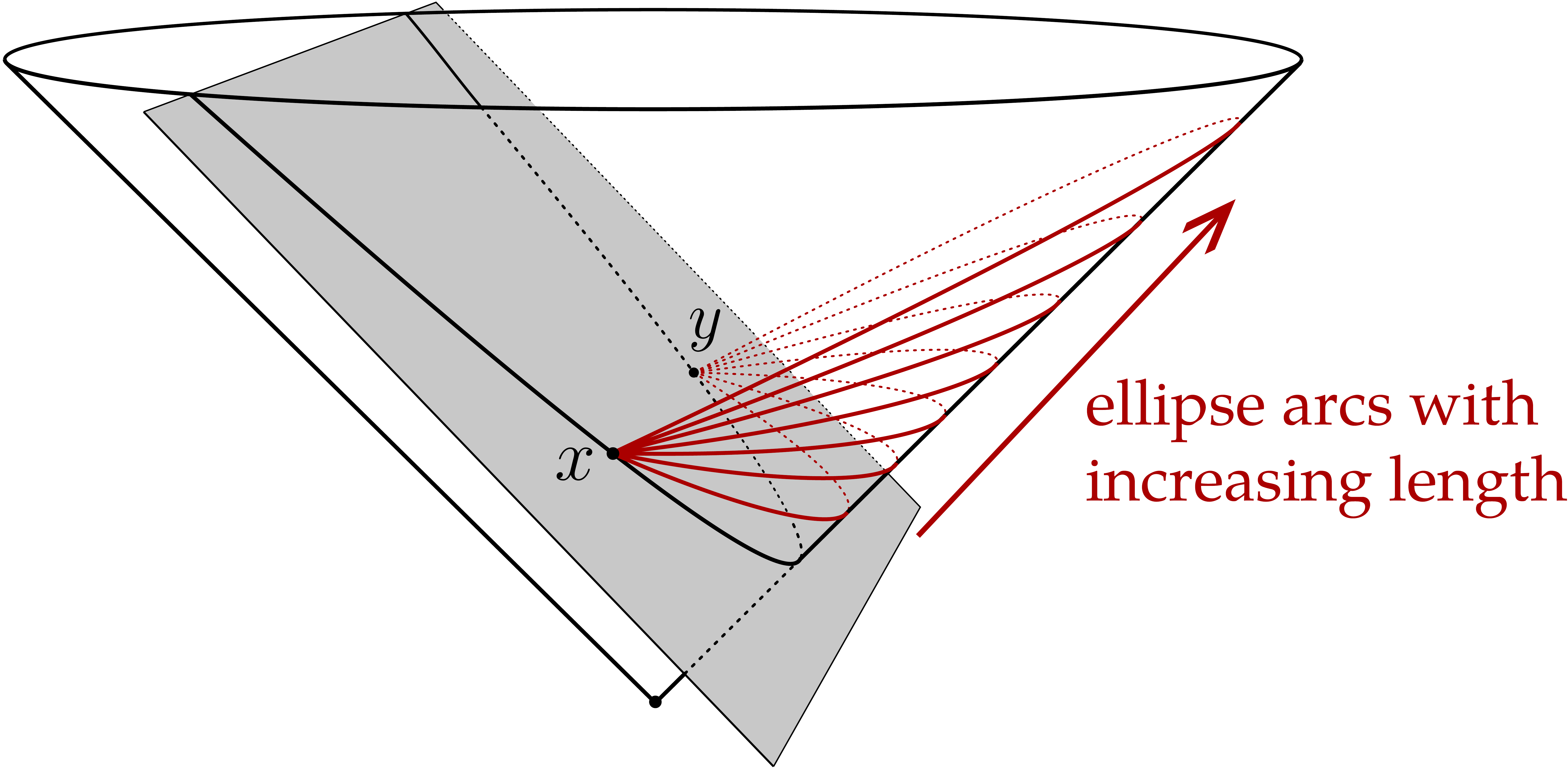}
	\caption{Property which implies ``\ref{item_Flemma1}$\Leftrightarrow$\ref{item_Flemma2}''.
	}
	\label{figure_arcs}
\end{figure}

In view of Lemma \ref{lemma_actionpair}, we may assume without loss of generality that 
\begin{equation}\label{eqn_standardxy}
x=(h,0,h),\quad y=(-h,0,h)
\end{equation}
for some $h>0$, and label the above arcs by the slope $k\in(-1,1)$ of the plane containing each of them (in Figure \ref{figure_arcs}, $k$ increases from $-1$ to $1$). Then the arc with slope $k$, denoted by $R_k$, is sent by the isometry
$$
\sigma_k:=
\begin{pmatrix}
\frac{1}{\sqrt{1-k^2}}&0&-\frac{k}{\sqrt{1-k^2}}\\
0&1&0\\
-\frac{k}{\sqrt{1-k^2}}&0&\frac{1}{\sqrt{1-k^2}}
\end{pmatrix}
\in\SO_0(2,1)
$$
to an arc $\sigma_k(R_k)$ on the horizontal circle $\big\{z\in\pa F\mid z_3=\frac{h}{\sqrt{1-k^2}}\big\}$. In fact, $\sigma_k(R_k)$ consists of points $z$ on this circle with second coordinate $z_2\geq-\frac{kh}{\sqrt{1-k^2}}$, which occupies a proportion of $\frac{\arccos(-k)}{\pi}$ in the full circle. By Lemma \ref{lemma_spacelike}, the full circle has perimeter $\frac{2\pi h}{\sqrt{1-k^2}}$. Therefore, we obtained the length of $R_k$ as
$$
\ell(R_k)=\ell(\sigma_k(R_k))=\frac{2h\arccos(-k)}{\sqrt{1-k^2}} 
$$
(the same under all three metrics), which is a strictly increasing function in $k\in(-1,1)$, as required.

\vspace{0.2cm}

\textbf{``\ref{item_Flemma2}$\Leftrightarrow$\ref{item_Flemma3}''.} Consider the point $o:=(0,0,1)$ on the hyperboloid $\H^2$. We note that for any $x\in\pa F\setminus\{0\}$, the signed distance in $\H^2$ from $o$ to the horocycle $x^*$ only depends on the vertical coordinate $x_3$ of $x$ and is strictly increasing in $x_3$. Indeed, by elementary calculations, this distance is $\log x_3$.

Given $\Delta_1''$ and $\Delta_2''$, by applying an element of $\SO_0(2,1)$ to both of them, we may assume that $\Delta_1''$ is circumscribed by a horizontal disk $\{x\in F\mid x_3=\text{constant}\}$. By the above noted fact, $o$ is the point with the same signed distance to all the horocycles of $\Delta_1'$, and \ref{item_Flemma3} can be reformulated as the condition that every vertex of $\Delta_2''$ lies no lower than that disk. This condition is clearly equivalent to \ref{item_Flemma2}, as required.

Finally, the ``Moreover'' statement is proved by following the same arguments, and we omit the details.
\end{proof}

\begin{lemma}\label{lemma_Ftrigonometry}
For any $x,y\in\pa F\setminus\{0\}$, the line segment $[x,y]\subset F$ is spacelike under the metrics $\ve{g}_{-1}$ and $\ve{g}_0$ of $\H^{2,1}$ and $\R^{2,1}$, and is spacelike under the metric $\ve{g}_{-1}$ of $\SS^{2,1}$ as well when $[x,y]\subset F^*$. If we let $\ell_{-1},\ell_0,\ell_1>0$ denote the lengths of $[x,y]$ under the three metrics respectively (the last is defined when $[x,y]\subset F^*$), and let
	$\ell'\in\R$ denote the signed distance between the horocycles $x^*$ and $y^*$ in $\H^2$, then we have
	\begin{equation}\label{eqn_Ftrigonometry}
	e^{\frac{\ell'}{2}}=\sqrt{-\frac{x\pro y}{2}}=
	\begin{cases}
		\sinh\frac{\ell_{-1}}{2},\\
		\tfrac{\ell_0}{2},\\
		\sin\frac{\ell_{1}}{2}.
	\end{cases}
	\end{equation}
\end{lemma}
\begin{proof}
 Again,	by Lemma \ref{lemma_linear}, we may assume  \autoref{eqn_standardxy} without loss of generality, where $h=\sqrt{-\frac{x\pro y}{2}}>0$. Then the spacelikeness follows immediately from the discussion at the end of \S \ref{subsec_preliminary}, while \autoref{eqn_Ftrigonometry} follows from Lemma \ref{lemma_spacelike} because $\frac{\ell_K}{2}$ is the radius of the horizontal disk $L_{(0,0,h)}\cap F$ under the metric $\ve{g}_K$ ($K=-1,0,1$). 
\end{proof}

Combining Lemmas \ref{lemma_trigonometry} and \ref{lemma_Ftrigonometry}, we get:
\begin{corollary}\label{coro_maincoro}
Given a convex cyclic polygon $\Delta$ in $\E^2$, $\H^2$ or $\SS^2$, let $\Delta''$ be a polygon in $F$ given by Lemma \ref{lemma_realizetriangle} (which isometrics to $\Delta$ when endowed with the metric induced from $\ve{g}_0$, $\ve{g}_{-1}$ or $\ve{g}_1$), and $\Delta'$ be the ideal polygonal face of the polyhedron $P(\Delta)\subset\H^3$, developed as an ideal polygon in $\H^2$ and endow with its canonical decoration (see the proof of Theorem \ref{thm_BPS}). Then $\Delta'$ is isomorphic to the decorated ideal polygon underlying $\Delta''$ in the cases of $\H^2$ and $\SS^2$. In the case of $\E^2$, $\Delta'$ is isomorphic to the decorated ideal polygon underlying $2\Delta'':=\big\{2x\,\big|\,x\in\Delta''\big\}$. 
\end{corollary}
The factor $2$ in the case of $\E^2$ comes from the fact that the equalities in Lemmas \ref{lemma_trigonometry} and \ref{lemma_Ftrigonometry} are the same for $\H^2$ and $\SS^2$ but differ by a factor of $2$ for $\E^2$. 
\begin{proof}
Label the sides of $\Delta$ by $1,2,\cdots,N$ in cyclic order and put the same labels on the corresponding sides of $\Delta''$ and $\Delta'$. For each $i\in\{1,\cdots,N\}$, let 
\begin{itemize}
\item 
$\ell_i$ be the length of the $i$th side of $\Delta$;
\item
$\ell'_i$ be the signed distance between the horocycles in $\H^2$ whose corresponding points in $\pa F\setminus\{0\}$ are the two ends of the $i$th side of $\Delta''$;
\item  
$\wt{\ell}'_i$ be the signed distance between the decorating horocycles at the two ends of the $i$th side of $\Delta'$.
\end{itemize}
By Lemmas \ref{lemma_Ftrigonometry} and \ref{lemma_trigonometry} respectively, $\ell'_i$ and $\wt{\ell}'_i$ are related to $\ell_i$ by 
$$
	e^{\frac{\ell'_i}{2}}=
\begin{cases}
	\sinh\frac{\ell_i}{2}\\
	\tfrac{\ell_i}{2}\\
	\sin\frac{\ell_i}{2}
\end{cases}
\quad\quad
	e^{\frac{\wt{\ell}'_i}{2}}=
\begin{cases}
	\sinh\frac{\ell_i}{2}\\
	\ell_i\\
	\sin\frac{\ell_i}{2}
\end{cases}
$$
where the three rows are for the cases $\H^2$, $\E^2$ and $\SS^2$, respectively. Therefore, we have $\ell'_i=\wt{\ell}'_i$ ($i=1,\cdots,N$) for $\H^2$ and $\SS^2$, while for $\E^2$ the same holds after scaling $\Delta''$ by a factor of $2$. The required statement follows because a convex cyclic polygon is determined by its side lengths (see e.g.\ \cite{Guo-Sonmez}).
\end{proof}

\begin{remark}\label{remark_chart2}
For $\H^{2,1}$ and $\SS^{2,1}$, our discussions still make sense when $F$ or $F^*$ is replaced by the future cone of another point and/or placed in another affine chart. In particular, if we relocate to the affine chart in Remark \ref{remark_chart1}, then $F$ gets transformed by the map \autoref{eqn_involution} to the semi-infinite cylinder $\wt{F}=\{x\in\R^3\mid x_1^2+x_2^2\leq1,\ x_3>0\}$, wherein the part $\wt{F}^*$ lying above the unit hemisphere corresponds to $F^*$. $\wt{F}$ and $\wt{F}^*$ are the future cone of a point at infinity for $\Omega_{-1}\subset\H^{2,1}$ and $\wt{\Omega}_1\subset\SS^{2,1}$, respectively.
A polygon $\Delta''\subset F$ as considered above corresponds to a polygon in $\wt{F}$ circumscribed by a boundary ellipse. See Figure \ref{figure_halfpipe}.
\begin{figure}[h]
	\centering
	\includegraphics[width=5.1cm]{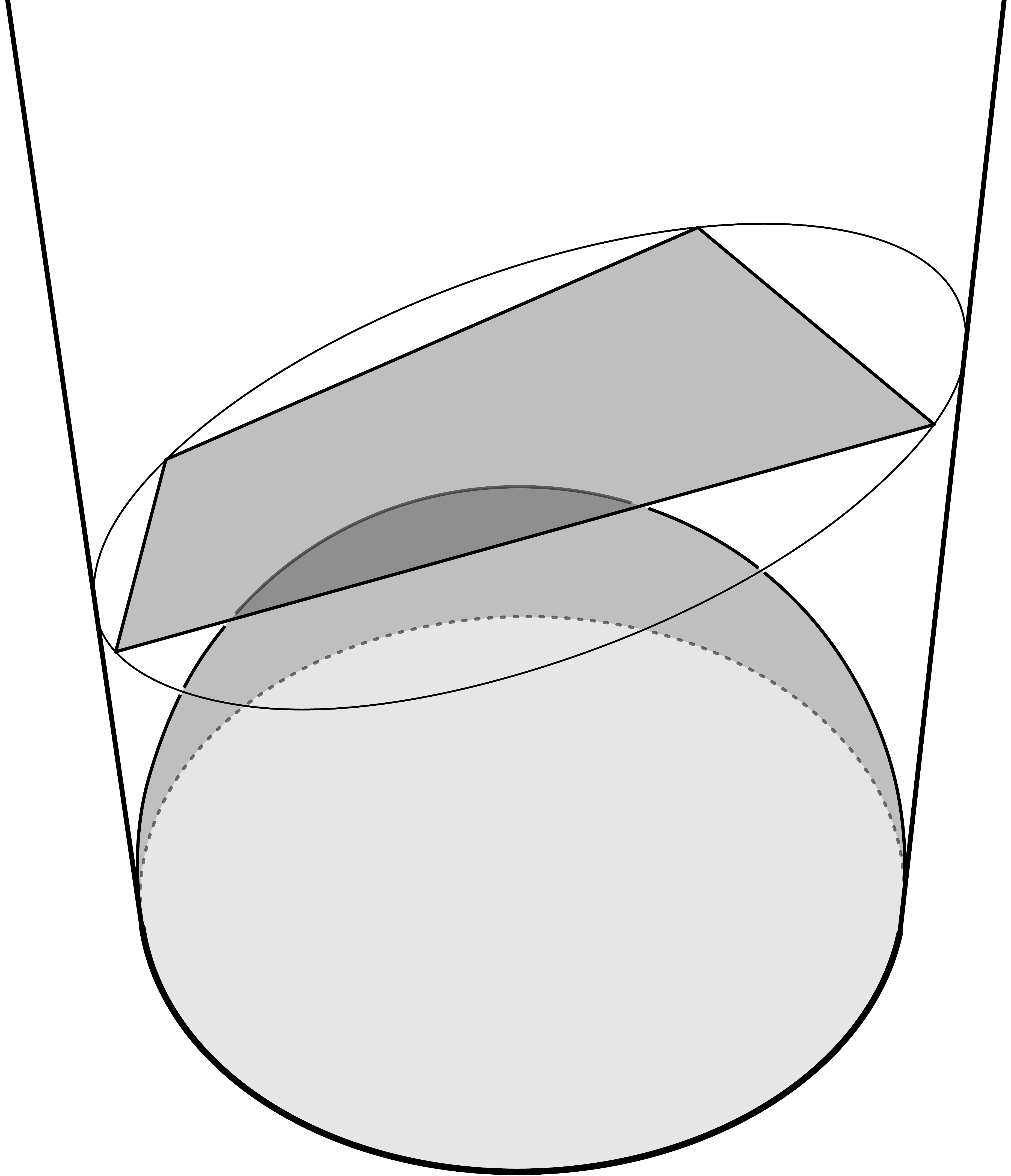}
	\hspace{2.3cm}
	\includegraphics[width=6.7cm]{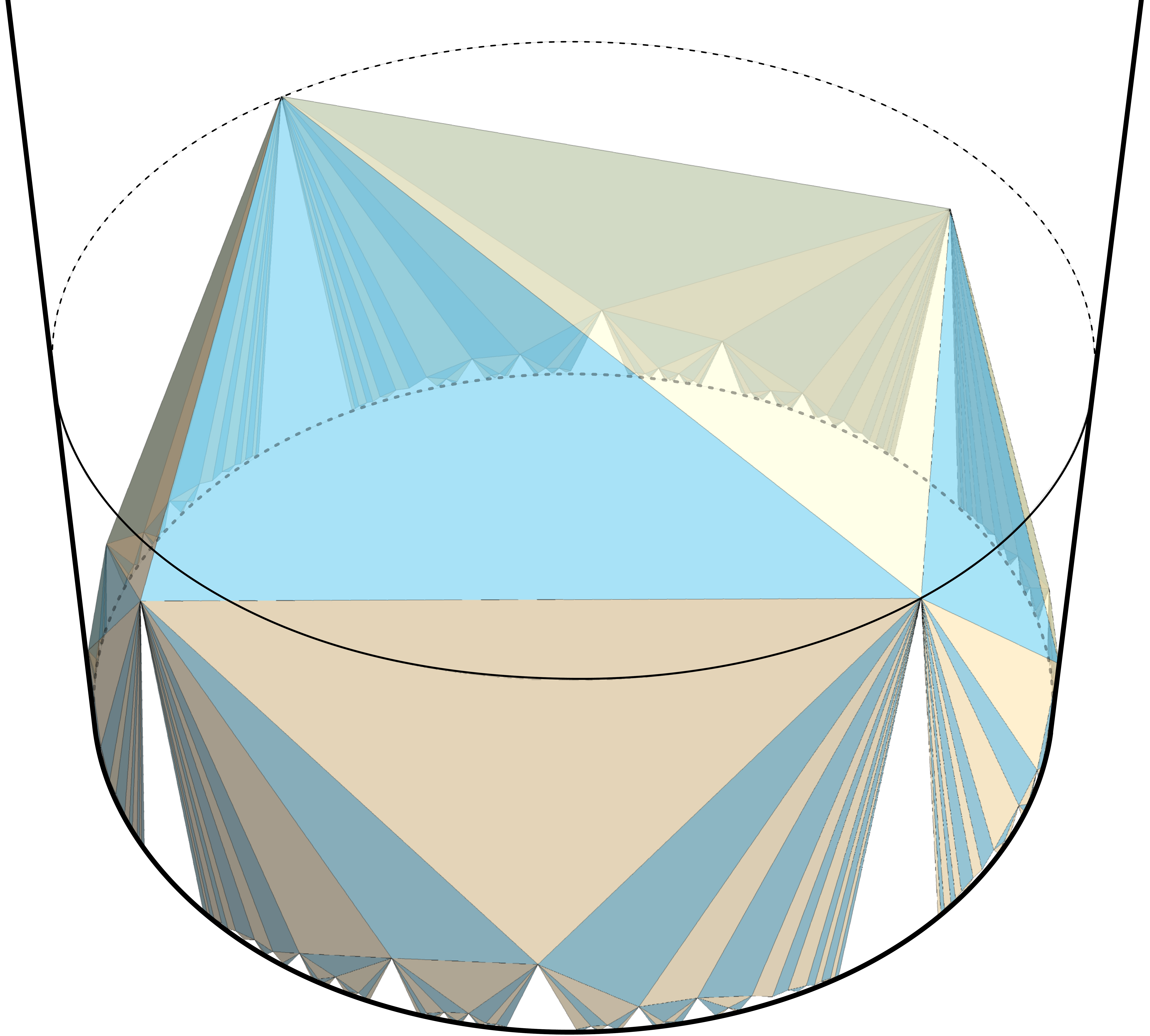}
	\caption{Left: future cone as a semi-infinite cylinder and a spacelike polygon in it. Right: an Epstein-Penner convex hull in the cylinder (transformed from the one in Figure \ref{figure_hull} by the projective transformation \autoref{eqn_involution}).
	}
	\label{figure_halfpipe}
\end{figure}
\end{remark}

\subsection{Epstein-Penner metrics}\label{subsec_EP}
Given a decorated hyperbolic metric $(\ve{h},c)\in\dT_{g,n}$, we have explained in Introduction the construction of the Epstein-Penner convex hull $\EP(\ve{h},c)\subset F$ (see Figures \ref{figure_hull} and \ref{figure_halfpipe}). Since $\EP(\ve{h},c)$ is well defined only up to the $\SO_0(2,1)$-action, a more rigorous notation is as follows. We identity $\TT_{g,n}$ as a subset in the space 
$\Hom(\pi_1(S_{g,n}^\circ),\SO_0(2,1))/\SO_0(2,1)$ 
formed by conjugacy classes of representations of the fundamental group $\pi_1(S_{g,n}^\circ)$ in $\SO_0(2,1)$. Then the convex hull is truly well defined after given the data of a point $(\ve{h},c)\in\dT_{g,n}$ along with a representation $\rho\in\Hom(\pi_1(S_{g,n}^\circ),\SO_0(2,1))$ which lifts $\ve{h}$, so we henceforth denote it instead by $\EP(\rho,c)$. Changing $\rho$ by some $A\in\SO_0(2,1)$ via conjugation amounts to changing $\EP(\rho,c)$ by the linear action of $A$.

The following properties of $\EP(\rho,c)$ are proved in Lemmas 3.3, 3.4 and Proposition 3.5 of \cite{Epstein-Penner}:
\begin{proposition}\label{prop_EPproperties}
Let $\rho\in\Hom(\pi_1(S_{g,n}^\circ),\SO_0(2,1))$ represent a point $\ve{h}$ of $\TT_{g,n}$, $c$ be a decoration of $\ve{h}$, and $C$ be the subset of $\pa F\setminus\{0\}$ corresponding to the lifts of the horocycles in $c$ (so that $\EP(\rho,c)$ is the convex hull of $C$). Then the boundary $\pa\EP(\rho,c)$ of $\EP(\rho,c)$ has the following properties:
\begin{enumerate}[label=(\arabic*)]
	\item\label{item_propEP1} Every ray in $F$ issuing from the origin but not lying on $\pa F$ meets $\pa E(\rho,c)$ exactly once.
	\item\label{item_propEP2} $\pa E(\rho,c)\cap \pa F$ is the union of rays $\{\alpha x\mid \alpha\geq1,\, x\in C\}$.
	\item\label{item_propEP3} $\pa^+\!E(\rho,c):=(\pa E(\rho,c)\setminus \pa F)\cup C$ is a union of planar polygons, where each polygon is circumscribed by an ellipse which is the intersection of $\pa F$ with an affine plane.
\end{enumerate}
\end{proposition}
These properties allow us to define the Euclidean (resp.\ hyperbolic) \emph{Epstein-Penner metric} $\ve{d}\in\M^0_{g,n}$ (resp.\ $\ve{d}\in\M^{-1}_{g,n}$) of $(\ve{h},c)$ as follows: By Part \ref{item_propEP1}, the radial projection maps $\pa E(\rho,c)\setminus \pa F$ diffeomorphically to the hyperboloid $\H^2$ in a $\rho$-equivariant way and descends to a diffeomorphism from the quotient of $\pa E(\rho,c)\setminus \pa F$ to the hyperbolic surface $(S_{g,n}^\circ,\ve{h})$. By \ref{item_propEP2}, if we add the points $C$ to the domain of the former map, then it descends to a diffeomorphism from the quotient of $\pa^+\!E(\rho,c)$ to the closed surface $S_{g,n}$. By \ref{item_propEP3} and Lemma \ref{lemma_spacelike}, the metric on $\pa^+\!E(\rho,c)$ induced by $\ve{g}_0$ (resp.\ $\ve{g}_{-1}$) is locally modeled on $\E^2$ (resp.\ $\H^2$), so the latter diffeomorphism turns it into the required metric $\ve{d}$ on $S_{g,n}$. Clearly, $\ve{d}$ depends only on $\ve{h}$ and not on the choice of its lift $\rho$. 

The spherical Epstein-Penner is defined in the same way, where we just need to be careful that ``the metric on $\pa^+\!E(\rho,c)$ induced by $\ve{g}_{1}$'' makes sense only when $\pa^+\!E(\rho,c)$ is contained in the region $\Omega_1$, or in other words, when $(\ve{h},c)$ is in the subset $\dT^*_{g,n}$ of $\dT_{g,n}$ defined in Introduction. Moreover, the resulting $\ve{d}\in\M^1_{g,n}$ belongs to the subspace $\M^{1,*}_{g,n}$ because if we subdivide the polygonal pieces of $\pa^+\!E(\rho,c)$ into triangles, then by Lemma \ref{lemma_spacelike}, each triangle has a circumdisk of radius less than $\frac{\pi}{2}$ in $\SS^2$, so they provide a triangulation of $\ve{d}$ into convex triangle. By Lemma \ref{lemma_FDelaunay}, this is actually a Delaunay triangulation, and the same is true in the Euclidean and hyperbolic cases.

Thus, we have obtained three Epstein-Penner metric assigning maps, which we will denote by
$$
\Psi_{-1}:\dT_{g,n}\to\M^{-1}_{g,n},\quad \Psi_0:\dT_{g,n}\to\M^{0}_{g,n},\quad \Psi_1:\dT^*_{g,n}\to\M^{1,*}_{g,n}.
$$
Given $(\ve{h},c)$ in $\dT_{g,n}$ or $\dT^*_{g,n}$, by the Delaunay property just mentioned, we have the following alternative, simpler description of the metric $\ve{d}=\Psi_K(\ve{h},c)$, which is similar in spirit to the definition of $\Phi_K$ in \S \ref{subsec_relationship polyhedral}: Letting $\T$ be the Delaunay decomposition or a Delaunay triangulation of $(\ve{h},c)$, we consider, for each face $\Delta'$ of $\T$, a polygon $\Delta''\subset F$ whose underlying decorated ideal polygon is $\Delta'$. Then $\ve{d}$ is obtained by converting every $\Delta'$ into the polygon $\Delta:=(\Delta'',\ve{g}_K|_{\Delta''})$, namely the hyperbolic/Euclidean/spherical polygon given by the metric on $\Delta''$ induced from $\ve{g}_K$.


This map $\Psi_K$ is essentially inverse to the map $\wt{\Phi}_K$ from Theorem \ref{thm_BPS}, modulo a dilation when $K=0$:
\begin{proposition}\label{prop_inverse}
For $K=\pm1$, $\Psi_K$ is inverse to $\wt{\Phi}_K$. Meanwhile, $\Psi_0$ is inverse to the composition $\wt{\Phi}_0\circ(\times \frac{1}{2})$, where ``$\times \frac{1}{2}$'' denotes the scaling map $\M^0_{g,n}\to\M^0_{g,n}$, $\ve{d}\mapsto \frac{1}{2}\ve{d}$.
\end{proposition}
\begin{proof}
In order to show, for example, that $\Psi_1\circ\wt{\Phi}_1$ is the identity self-map of $\M^{1,*}_{g,n}$, we pick $\ve{d}\in\M^{1,*}_{g,n}$, let $\T$ be its Delaunay decomposition, and recall from the proof of Theorem \ref{thm_BPS} that $(\ve{h},c)=\wt{\Phi}_1(\ve{d})$ is defined by converting every face $\Delta$ of $\T$ (a convex cyclic spherical polygon) into the face $\Delta'$ of $P(\Delta)$ (a decorated ideal polygon).
By Lemma \ref{lemma_Delaunay}, these $\Delta'$ also form a Delaunay decomposition of $(\ve{h},c)$. Therefore, as mentioned preceding this proposition, $\wt{\ve{d}}:=\Psi_1(\ve{h},c)=\Psi_1\circ\wt{\Phi}_1(\ve{d})$ is obtained by considering, for every $\Delta'$, the polygon $\Delta''\subset F^*$ whose underlying decorated ideal polygon is $\Delta'$, and then converting $\Delta'$ into the spherical polygon $\wt{\Delta}:=(\Delta'',\ve{g}_1|_{\Delta''})$. By Corollary \ref{coro_maincoro}, $\wt{\Delta}$ is isometric to $\Delta$.
So we have $\wt{\ve{d}}=\ve{d}$, which means that $\Psi_1\circ\wt{\Phi}_1=\id$, as required.
With a similar argument, one can show that $\wt{\Phi}_1\circ\Psi_1$ is also the identity self-map of $\dT^*_{g,n}$. Thus, we conclude that $\Psi_1$ is inverse to $\wt{\Phi}_1$. The cases $K=-1,0$ are similar.
\end{proof}
\begin{remark}\label{remark_fiberconnected}
Although not needed in this paper, it is easy to see that the fiber of $\dT^*_{g,n}$ over any $\ve{h}\in\TT_{g,n}$ is connected.
In fact, if we fix an auxiliary decoration $c=(c_1,\cdots,c_n)$ of $\ve{h}$
and identify the fiber $\Pi^{-1}(\ve{h})$ in $\dT_{g,n}$ with $\R^n$ in such a way that $s=(s_1,\cdots,s_n)\in\R^n$ corresponds to the decoration obtained by shrinking each horocycle $c_i$ by a signed distance of $s_i$, or equivalently, by scaling each point in $\pa F\setminus\{0\}$ corresponding to a lift of $c_i$ by a factor of $e^{s_i}$, then $\dT^*_{g,n}\cap\Pi^{-1}(\ve{h})$ corresponds to a subset $S\subset\R^n$ with the property (stronger than connectedness) that for any $s\in S$, every point $\wh{s}\in\R^n$ coordinate-wise less than or equal to $s$ is again in $S$. This is simply because we have $E(\rho,\wh{s})\subset E(\rho,s)$ for any $\rho$ lifting $\ve{h}$.    
\end{remark}

\subsection{Proof of Theorem \ref{thm_intro}}
With the map $\Psi_K$ introduced above, the theorem can be restated as:
\begin{theorem}
For each $K\in\{-1,0,1\}$, the map $\Psi_K$ is bijective and sends each fiber of $\Pi$ in its domain bijectively to a discrete conformal class in its target.
\end{theorem}
This can largely be deduced from Theorem \ref{thm_BPS} and Proposition \ref{prop_inverse}. For the sake of completeness, we give here a direct proof. Much of the argument is similar to that in the proof of  Theorem \ref{thm_BPS}.

\begin{proof}
\textbf{Surjectivity.} We shall show that any $\ve{d}\in\M^{-1}_{g,n}$, $\M^0_{g,n}$ or $\M^{1,*}_{g,n}$ is the Epstein-Penner metric of some $(\ve{h},c)\in\dT_{g,n}$. To this end, we let $\T$ be the Delaunay decomposition of $\ve{d}$ and construct $(\ve{h},c)$ as follows, where $\ve{h}$ will be given by a representation $\rho:\pi_1(S_{g,n}^\circ)\to \SO_0(2,1)$ as explained at the beginning of \S \ref{subsec_EP}.

Consider the lifts $\wt{\T}$ and $\wt{\ve{d}}$ of $\T$ and $\ve{d}$ to the universal cover of $S^\circ_{g,n}$. For each face $\Delta$ of $\wt{\T}$, we pick a spacelike planar polygon $\Delta''$ in $F$ isometric to $\Delta$ as provided by Lemma \ref{lemma_realizetriangle}. By using Lemma \ref{lemma_actionpair}, we may assume that these $\Delta''$$\,$'s are pieced together in the same combinatorics as the $\Delta$$\,$'s to form a piecewise linear surface $\Sigma\subset F$. Those two lemmas also imply that once a single $\Delta''$ is chosen, all the others are uniquely determined, and furthermore there exists a representation $\rho:\pi_1(S^\circ_{g,n})\to \SO_0(2,1)$ such that if two faces $\Delta_1$ and $\Delta_2$ of $\wt{\T}$ are related by a deck transformation $\gamma\in\pi_1(S^\circ_{g,n})$, then $\Delta''_1$ and $\Delta''_2$ are related by $\rho(\gamma)\in\SO_0(2,1)$. Moreover, since $\T$ is Delaunay for $\ve{d}$, the surface $\Sigma$ is convex by Lemma \ref{lemma_FDelaunay}. Therefore, $\rho$ and the vertices of the $\Delta''$$\,$'s define respectively a point $\ve{h}$ in $\dT_{g,n}$ and a decoration $c$ of $\ve{h}$, such that $\Sigma$ is exactly the spacelike boundary part $\pa^+\!E(\rho,c)$ of the Epstein-Penner convex hull $\EP(\rho,c)$. It follows that $\ve{d}$ is the Epstein-Penner metric of $(\ve{h},c)$, as required.

\vspace{0.1cm}

\textbf{Injectivity.} We shall show that if $\ve{d}\in\M^{-1}_{g,n}$, $\M^0_{g,n}$ or $\M^{1,*}_{g,n}$ is the Epstein-Penner metric of some $(\ve{h},c)\in\dT_{g,n}$, then $(\ve{h},c)$ is uniquely determined by $\ve{d}$. 
By Lemma \ref{lemma_FDelaunay}, the Delaunay decomposition of $\ve{d}$ is also the Delaunay decomposition of $(\ve{h},c)$, so $\ve{d}$ determines the combinatorics of the Delaunay decomposition of $(\ve{h},c)$. By Lemma \ref{lemma_Ftrigonometry}, the edge lengths of the Delaunay decomposition of $\ve{d}$ also determine the Penner coordinates of $(\ve{h},c)$ at the edges of its Delaunay decomposition. But every point $(\ve{h},c)$ in $\dT_{g,n}$ is determined by the combinatorics of its Delaunay decomposition $\T$ together with the Penner coordinate $\ell'_e$ at all $e\in\Edge_\T$, because any cyclic decorated ideal polygon is determined by the signed distances between the horocycles at all adjacent vertices. Therefore, we conclude that $\ve{d}$ determines $(\ve{h},c)$, as required.

\vspace{0.1cm}

\textbf{Characterization of discrete conformality.}
As mentioned in the proof of Theorem \ref{thm_BPS}, two points of $\dT_{g,n}$ are in the same fiber of $\Pi$ if and only if they are respectively the first and last members of a finite sequence in $\dT_{g,n}$ such that any adjacent members $(\ve{h},c)$ and $(\wh{\ve{h}},\wh{c})$ of the sequence
satisfy the following condition for some topological triangulation $\T$ of $S_{g,n}$:
\begin{equation}\label{eqn_proofmain3}\tag{$\star'_{\scalebox{0.5}{$\T$}}$}
	\parbox{12.5cm}{
		$(\ve{h},c)$ and $(\wh{\ve{h}},\wh{c})$ are both in the intersection of the Penner cell $\dT^\T_{g,n}$ with some fiber of $\Pi$; or equivalently, they are both in $\dT_{g,n}^\T$ and there exists $(u_i)\in\R^n$ such that their Penner coordinates $(\ell'_e)_{e\in\Edge_\T}$ and $(\wh{\ell}'_e)_{e\in\Edge_\T}$ are related by $\hat{\ell}'_e=\ell'_e+u_i+u_j$ for all $i,j\in\{1,\cdots,n\}$ and all $e\in\Edge_\T$ joining the $i$th and $j$th marked points.
	}
\end{equation}
Meanwhile, two metrics in $\M^{-1}_{g,n}$, $\M^0_{g,n}$ or $\M^{1,*}_{g,n}$ are discretely conformal if and only if they are the first and last members of a finite sequence wherein any adjacent members $\ve{d}$ and $\wh{\ve{d}}$ 
satisfy the following for some $\T$: 
\begin{equation}\label{eqn_proofmain1}\tag{$\star_{\scalebox{0.5}{$\T$}}$}
	\parbox{13.2cm}{
		$\T$ can be realized as Delaunay triangulations for both $\ve{d}$ and $\wh{\ve{d}}$. Moreover, there exists $(u_i)\in\R^n$ such that 
		for all $i,j\in\{1,\cdots,n\}$ and all $e\in\Edge_\T$ joining the $i$th and $j$th marked points, the geodesic lengths $\ell_e$ and $\hat{\ell}_e$ of $e$ under $\ve{d}$ and $\wh{\ve{d}}$ are related by  
		$$
		\hspace{-3cm}
		\begin{cases}
			\sinh\tfrac{\hat\ell_e}{2}=e^{\frac{u_i+u_j}{2}}\sinh\tfrac{\ell_e}{2}&\text{ for $\M^{-1}_{g,n}$,}\\[5pt]    
			\hat\ell_e=e^{\frac{u_i+u_j}{2}}\ell_e&\text{ for $\M^{0}_{g,n}$,}\\[5pt]
			\sin\tfrac{\hat\ell_e}{2}=e^{\frac{u_i+u_j}{2}}\sin\tfrac{\ell_e}{2}&\text{ for $\M^{1,*}_{g,n}$.}
		\end{cases}
		$$
	}
\end{equation}
By the above mentioned coincidence between the Delaunay decomposition of any $(\ve{h},c)\in\dT_{g,n}$ and that of its image $\ve{d}=\Psi_K(\ve{h},c)$ which follows from Lemma \ref{lemma_FDelaunay}, as well as the relation between Penner coordinates of $(\ve{h},c)$ and the edge lengths of $\ve{d}$ which follows from Lemma \ref{lemma_Ftrigonometry}, we infer that $(\ve{h},c)$ and $(\wh{\ve{h}},\wh{c})$ satisfy \autoref{eqn_proofmain3} if and only if $\ve{d}=\wt{\Phi}_K(\ve{h},c)$ and $\wh{\ve{d}}=\wt{\Phi}_K(\wh{\ve{h}},\wh{c})$ satisfy \autoref{eqn_proofmain1}. It follows that any  $(\ve{h}_1,c_1),(\ve{h}_2,c_2)\in\dT_{g,n}$ are in the same fiber of $\Pi$ if and only if $\ve{d}_1=\Psi_K(\ve{h}_1,c_1)$ and $\ve{d}_2=\Psi_K(\ve{h}_2,c_2)$ are discretely conformal, as required.
\end{proof}

\end{document}